\documentclass[leqno,11pt]{article}
\usepackage{amssymb,latexsym}
\usepackage[dvips]{graphicx} 
\usepackage{url}

 \catcode`@=11\def\@oddhead{\hbox{}\hfil\rm\thepage}\def\@oddfoot{}
 \def\@evenhead{\hbox{}\hfil\rm\thepage}\def\@evenfoot{}
 \@twosidetrue \catcode`@=12
\newtheorem{prp}{Proposition}
\newtheorem{lem}[prp]{Lemma}\newtheorem{thm}[prp]{Theorem}

\newenvironment{prf}{\begin{trivlist}\item[\emph{Proof.}]}{\end{trivlist}
  \medskip\par}
\newenvironment{prfof}[1]{\begin{trivlist}\item[\emph{Proof of #1.}]}{
  \end{trivlist} \medskip \par}
\newenvironment{rem}{\begin{trivlist}\item[\emph{Remark.}]}{\end{trivlist}
  \medskip\par}
\def\prpa#1{\label{p:#1}}\def\prpu#1{Proposition~\ref{p:#1}}
\def\lema#1{\label{l:#1}}\def\lemu#1{Lemma~\ref{l:#1}}
\def\thma#1{\label{t:#1}}\def\thmu#1{Theorem~\ref{t:#1}}

\def\seca#1{\label{s:#1}}\def\secu#1{\S~\ref{s:#1}}

\def\eqna#1{\label{e:#1}}\def\eqnu#1{(\ref{e:#1})}

\def\QED{\relax\ifmmode\let\@tempa\relax\ifcase\@eqcnt\def\@tempa{& & &}\or
  \def\@tempa{& &}\else\def\@tempa{&}\fi\@tempa $\Box$ \else\hfill $\Box$ \fi}
\def\DDD{\relax\ifmmode\let\@tempa\relax\ifcase\@eqcnt\def\@tempa{& & &}\or
 \def\@tempa{& &}\else\def\@tempa{&}\fi\@tempa $\Diamond$
 \else\hfill $\Diamond$ \fi}

\def\Rom#1{\uppercase\expandafter{\romannumeral#1}}
\def\dsp{\displaystyle}
\def\eps{\epsilon}

\def\limf#1{\displaystyle \lim_{#1\to\infty}}

\def\limsup{\displaystyle \mathop{\overline{\lim}}\limits}

\def\diff#1#2{\dsp\frac{d\,#1}{d#2}}

\def\pderiv#1#2{\dsp\frac{\partial\,#1}{\partial#2}}

\def\le{\leqq} \def\ge{\geqq} 

\def\reals{{\mathbb R}}

\def\preals{[0,\infty)} 
\def\pintegers{{\mathbb Z}_+}

\def\prb#1{\def\prbone{#1}
  \ifx\prbone\empty{\mathrm{P}}\else{\mathrm{P[\;}}#1{\mathrm{\;]}}\fi}
\def\prbseq#1#2{\def\prbseqone{#2}
  \ifx\prbseqone\empty{\mathrm{P}}_{#1}\ignorespaces
  \else{\mathrm{P}}_{#1}{\mathrm{[\;}}#2{\mathrm{\;]}}\fi}
\def\EE#1{{\mathrm{E[\;}}#1{\mathrm{\;]}}}
\def\EEseq#1#2{\def\EEseqone{#2}
  \ifx\EEseqone\empty{\mathrm{E}}_{#1}\else
 {\mathrm{E}}_{#1}{\dsp\mathrm{[\;}}#2{\mathrm{\;]}}\fi}

\def\VVseq#1#2{\def\VVseqone{#2}
  \ifx\VVseqone\empty{\matrm{V}}_{#1}\else
 {\mathrm{V}}_{#1}{\dsp\mathrm{[\;}}#2{\mathrm{\;]}}\fi}

\def\ss#1{^{(#1)}}

\def\chrfcn#1{\mathop{\mathbf{1}}\nolimits_{#1}}

\def\wnorm#1{\dsp\left\| #1 \right\|_{\rm{T}}}

\title{
Point process with last-arrival-time dependent intensity
and $1$-dimensional incompressible fluid system with evaporation.
}
\author{
Tetsuya Hattori
\footnote{ 
This work is supported by JSPS KAKENHI Grant Number 26400146.
}
\\ 
\small Laboratory of Mathematics, Faculty of Economics, Keio University, 
\\
\small Hiyoshi Campus, 4--1--1 Hiyoshi, Yokohama 223-8521, Japan
\\ \small URL: \url{http://web.econ.keio.ac.jp/staff/hattori/research.htm}
\\ \small email: \url{hattori@econ.keio.ac.jp}
} 
\date{\today}

\begin{document}

\setcounter{section}{0}

\setcounter{equation}{0}

\maketitle

\begin{abstract}
We consider an infinite system of quasilinear
first-order partial differential equations, generalized to contain
spacial integration, which describes an incompressible 
fluid mixture of infinite components in a line segment
whose motion is driven by unbounded and space-time dependent
evaporation rates.
We prove unique existence of the solution to the initial-boundary value
problem, with conservation-of-fluid condition at the boundary.
The proof uses a map on the space of collection of characteristics,
and a representation based on a non-Markovian point process
with last-arrival-time dependent intensity.
\end{abstract}

MSC2010-class: Primary~35F61; Secondary~45G15, 35C99.

Keywords:
quasilinear first-order infinite system, non-Markov point process,
partial differential integral equations

\section{Introduction.}
\seca{intro}

Consider an incompressible fluid mixture in a line segment,
say $[0,1]$, which flow in order preserving manner and in one direction,
with $y=0$ being the upper stream boundary,
and no leaking occurs at $y=1$.
Each fluid component, say $\alpha$, evaporates with rate $w_{\alpha}$
which may vary among different components
and may depend on time.
Flow of the fluid is driven by filling the evaporated portion
of the fluid toward the down stream.
To formulate a system of partial differential equations
which explains the dynamics of this fluid up to time $T>0$,
let $U_{\alpha}(y,t)$ be the total volume (length) of fluid component
$\alpha$ at time $t$ in the interval $[y,1)$.
Then we have
\begin{equation}
\eqna{Burgersspaceindep}
\pderiv{U_{\alpha}}{t}(y,t) 
+\sum_{\beta}w_{\beta}(t)U_{\beta}(y,t)\, \pderiv{U_{\alpha}}{y}(y,t)
=-w_{\alpha}(t) U_{\alpha}(y,t),\ \ (y,t)\in[0,1]\times[0,T].
\end{equation}
We will preserve the total volume of each component by supplying
the evaporated portion from upper stream boundary through the
boundary condition
\begin{equation}
\eqna{Burgersconservation}
U_{\alpha}(0,t)=r_{\alpha}\ \mbox{ and }\ U_{\alpha}(1,t)=0,\ \ t\ge0,
\end{equation}
for non-negative constants $r_{\alpha}$ satisfying 
$\dsp \sum_{\alpha}r_{\alpha}=1$.
The incompressibility condition is formulated as
\begin{equation}
\eqna{Burgerssoliditycond}
\sum_{\alpha} U_{\alpha}(y,t)=1-y,\ \ t\ge0.
\end{equation}
The number of fluid components may be finite or infinite.
(For the latter case we regard the summations in $\alpha$ as series.)
If the system is infinite, we should impose an additional condition
\begin{equation}
\eqna{Burgersspaceindeprw}
\sum_{\alpha} r_{\alpha} \sup_t w_{\alpha}(t)<\infty,
\end{equation}
to keep the velocity of the flow finite,
namely, to keep coefficient of the $y$-derivative term in
\eqnu{Burgersspaceindep} well-defined at $y=0$.
With appropriate initial conditions, these define 
an initial/boundary value problem of a one dimensional 
first order quasilinear partial differential equations.

In this paper we consider a generalization of \eqnu{Burgersspaceindep}
to allow for spacial dependence for the evaporation rates $w_{\alpha}$,
as well as time dependence.
Such generalization seems practically natural,
because if the fluid container has spacial non-uniformity in temperature, 
the evaporation rates would also have spacial dependence.
$U_{\alpha}(y,t)$ is the volume of type $\alpha$ fluid component
in the interval $[y,1)$, and we need to consider its density
to consider spacially varying evaporation rates,
hence a natural generalization of \eqnu{Burgersspaceindep} is
\begin{equation}
\eqna{Burgers}
\begin{array}{l}\dsp
\pderiv{U_{\alpha}}{t}(y,t) 
-\sum_{\beta}\int_y^1 w_{\beta}(z,t)\,\pderiv{U_{\beta}}{z}(z,t)\,dz\,
\pderiv{U_{\alpha}}{y}(y,t)
=\int_y^1 w_{\alpha}(z,t)\,\pderiv{U_{\alpha}}{z}(z,t)\,dz,
\\ \dsp
\alpha=1,2,\ldots,\ (y,t)\in[0,1]\times[0,T].
\end{array}
\end{equation}
Note that the equation is now non-local and contains integration.
If $w_{\alpha}$ are independent of $y$, then, 
with \eqnu{Burgersconservation}, \eqnu{Burgers} reduces to
\eqnu{Burgersspaceindep}.

The equation of the form \eqnu{Burgersspaceindep} is known to
be solved by considering characteristic curves \cite{Bressan},
a curve $y=y_C(t)$ whose derivative is equal to the velocity 
of fluid;
\begin{equation}
\eqna{Burgersspaceindepfluidvelocity}
\diff{y_C}{t}(t)=\sum_{\beta} w_{\beta}(t)\,\varphi_{\beta}(t)
\end{equation}
where $\dsp\varphi_{\alpha}(t)=U_{\alpha}(y_C(t),t)$.
Then \eqnu{Burgersspaceindep} implies 
an ordinary differential equation for $\varphi_{\alpha}(t)$,
which can be solved explicitly, and
\eqnu{Burgersspaceindepfluidvelocity} then implies
\begin{equation}
\eqna{BurgersspaceindepyC}
y_C(t)=1-\sum_{\beta}U_{\beta}(y_0,0)\, \exp(-\int_{t_0}^t w_{\beta}(u)\,du).
\end{equation}
A natural generalization of \eqnu{Burgersspaceindepfluidvelocity} for
\eqnu{Burgers} is
\begin{equation}
\eqna{fluidvelocity}
\diff{y_C}{t}(t)=
-\sum_{\beta}\int_{y_C(t)}^1
w_{\beta}(z,t)\,\pderiv{U_{\beta}}{z}(z,t)\,dz.
\end{equation}
As we will see in the present paper,
this is no longer solved in such simple form as \eqnu{BurgersspaceindepyC}.
Introduction of spacial dependence for $w_{\alpha}$ complicates
the solution when combined with the boundary condition 
\eqnu{Burgersconservation} which conserves component volumes.
We will show later that $y_C$ is determined as a fixed point 
to the map $G$ defined by \eqnu{HDLop_on_class_of_flows_f} and
\eqnu{HDLop_on_class_of_flows_g}, a result which apparently deviates
largely from \eqnu{BurgersspaceindepyC}.
To be specific, our proof in \secu{FP} of \thmu{FPtheorem} 
proves an expression $\dsp y_C=\limf{n} G^n(\theta_0)$
for the characteristic curves, where $\theta_0$ is
a constant flow.

In the preceding work \cite{posdep}, 
the problem \eqnu{Burgers} for $w_{\alpha}$ with spacial dependence
was considered under the condition
\[
\sup_{\alpha}\sup_{(y,t)} w_{\alpha}(y,t)<\infty,
\ \mbox{ and }\ 
\sup_{\alpha}\sup_{(y,t)} \pderiv{w_{\alpha}}{y}(y,t)<\infty.
\]
In view of \eqnu{Burgersspaceindeprw} for $\{w_{\alpha}\}$ with spacially
independent case, a natural restriction for $w_{\alpha}$
is expected to be a milder one,
\begin{equation}
\eqna{Burgersrw}
\sum_{\alpha} r_{\alpha} \sup_{(y,t)} w_{\alpha}(y,t)<\infty,
\ \mbox{ and }\ 
\sup_{\alpha}\sup_{(y,t)} \pderiv{w_{\alpha}}{y}(y,t)<\infty,
\end{equation}
allowing, in particular, fluid mixture with unbounded evaporation rates.
The unique existence of the solution was proved for the case of
bounded evaporation rates in \cite{posdep},
but explicit formula such as \eqnu{HDLop_on_class_of_flows_g} were absent.
A main interest in \cite{posdep} is on the stochastic ranking process 
\cite{HH071,HH072,timedep,Nagahata10,Nagahata13},
and a rather strong restriction on $w_{\alpha}$ was posed to prove
existence of the hydrodynamic limit of the process
with spacially dependent $w_{\alpha}$, and \eqnu{Burgers} appears as
the equation which characterizes the limit of the process.
In the present paper we focus on \eqnu{Burgers} itself,
and solve the equation under a natural assumption \eqnu{Burgersrw}.
Besides mathematical naturalness of the assumption \eqnu{Burgersrw},
removal of boundedness condition on $\{w_{\alpha}\}$ has practical
meaning also on an application of stochastic ranking process
to an analysis on a behavior of web ranking data for on-line retail 
businesses.
See \cite{ranking,mv2frnt,dojin}, in addition to the references above
for details on practical applications to web ranking.

As we will show in this paper,
the solution $U_{\alpha}$ turns out to have a concise expression
using the stochastic processes $N_{\theta,w,z}$
which we introduce in \secu{pointprocess},
\begin{equation}
\eqna{BurgersspaceindepU}
U_{\alpha}(y,t)= \int_{z\in[y_0,1)}
\prb{N_{y_C,w_{\alpha},z}(t)=N_{y_C,w_{\alpha},z}(t_0)}\,\mu_{0,\alpha}(dz),
\end{equation}
where $\mu_{0,\alpha}$ denotes initial spacial distribution of the fluid
component $\alpha$, and $(y_0,t_0)$ is a initial/boundary point such that
the characteristic curve starting from the point satisfies $y=y_C(t)$.
(See \eqnu{explicit_soln} with \eqnu{mu2U}.)
The map $G$ of \eqnu{HDLop_on_class_of_flows_f} and
\eqnu{HDLop_on_class_of_flows_g} also has a corresponding expression \eqnu{G}.
The processes $N_{\theta,w,z}$
may be regarded as generalizations of the Poisson process,
but, in contrast to the Poisson process,
lacks independent increment properties, 
resulting in the complexity of the solution.
In the case of spacially independent evaporation rates,
this underlying process reduces to the Poisson process,
whose independent increment property implies simple explicit formula
such as \eqnu{BurgersspaceindepyC}.

We mentioned earlier that the characteristic curve,
which is the key quantity for a solution to a
one dimensional first order quasilinear partial differential equation,
will no longer be obtained by ordinary differential equation
for the spacially dependent $\{w_{\alpha}\}$,
and that it is determined as the fixed point to a map.
The map is 
on the collection of the characteristic curves parametrized by
its intersection point with the initial/boundary points,
the totality of which we introduce as flow in \secu{flows}.

These notions were absent in the preceding work \cite{posdep},
and it is to clarify such mathematical structure of the solution
that mainly motivated the present paper.

The plan of the paper is as follows.
In \secu{theorem} we give the precise statement of our result,
where we generalize \eqnu{Burgers} to allow also for 
uncountable number of fluid components, by generalizing
the unknown functions to measure valued function.
In \secu{pointprocess} we introduce the underlying stochastic process
and its elementary properties, with which we give an expression 
of the solution in \secu{flows} (see \eqnu{explicit_soln}),
assuming existence of a fixed point to a certain map (\thmu{soln}).
The existence of the fixed point is proved in \secu{FP},
which completes the existence proof of the solution.
A uniqueness proof of the solution is given in \secu{prfofBurgers}.
As a remark concerning the condition in \eqnu{Burgersrw}
on the spacial derivatives of $w_{\alpha}$,
we apply Schauder's fixed point theorem in \secu{Schauder}
to the map defined by \eqnu{HDLop_on_class_of_flows_f} and 
\eqnu{HDLop_on_class_of_flows_g}, with the condition on derivative
relaxed to a global bound on oscillation of $w_{\alpha}$\,.

\section{Main Result.}
\seca{theorem}

Throughout this paper we fix $T>0$, 
$W\subset C^1([0,1]\times[0,T];\preals)$
a set of non-negative valued $C^1$ functions on $[0,1]\times[0,T]$,
and a Borel probability measure $\lambda$
supported on the Borel measurable space $(W,{\cal B}(W))$.
${\cal B}(W)$ is the $\sigma$-algebra generated by open sets
with the topology from the space of continuous functions
$C^0([0,1]\times[0,T];\preals)\supset C^1([0,1]\times[0,T];\preals)$
with the metric given by the supremum norm
\begin{equation}
\eqna{Wnorm}
\wnorm{w}=\sup_{(y,t)\in[0,1]\times[0,T]} |w(y,t)|.
\end{equation}
We assume that
\begin{equation}
\eqna{rw}
M_W:=\int_W \wnorm{w} \lambda(dw) <\infty
\end{equation}
and
\begin{equation}
\eqna{Cw}
C_W:= \sup_{w\in W} \wnorm{\pderiv{w}{y}} <\infty
\end{equation}
hold.

Denote the sets of `initial ($t=0$) points' 
in the space-time $[0,1]\times[0,T]$,
the set of `upper stream boundary ($y=0$) points',
and their union, the set of initial/boundary points, respectively by
\begin{equation}
\eqna{Gamma}
\begin{array}{l}\dsp
\Gamma_b=\{0\}\times [0,T]=\{(0,s)\mid 0\le s\le T\},
\\
\Gamma_i=[0,1]\times \{0\}=\{(z,0)\mid 0\le z\le 1\},
\\
\Gamma=\Gamma_b\cup\Gamma_i\,.
\end{array}
\end{equation}
For $t\in[0,T]$, denote 
the set of initial/boundary points up to time $t$ by
\begin{equation}
\eqna{Gammat}
\Gamma_t=\{(z,s)\in \Gamma \mid t_0\le t\}
=\Gamma_i\cup \{(0,t_0)\in\Gamma_b \mid 0\le t_0\le t\},
\end{equation}
and the set of admissible pairs of the initial/boundary point $\gamma$
and time $t$ by
\begin{equation}
\eqna{domyC}
\Delta_T:=\{(\gamma,t)\in\Gamma_T\times [0,T]\mid \gamma\in\Gamma_t\}.
\end{equation}

To state the initial condition,
let $\mu_0=\mu_0(dw\times dz)$ be a Borel probability measure
on the measurable space $(W\times[0,1],{\cal B}(W\times[0,1]))$
of the product space of $W$ and $[0,1]$.
We assume that $\mu_0$ is absolutely continuous with respect to
the product measure $\lambda\times dz$, where $dz$ denotes the
standard Lebesgue measure on $\reals$.
Denote the density function by $\sigma$, namely,
\begin{equation}
\eqna{Burgersdefdensity}
\mu_0(dw\times dz)=\sigma(w,z)\,\lambda(dw)\, dz,\ (w,z)\in W\times[0,1].
\end{equation}
We assume $\mu_0(W\times dz)=dz$ and $\mu_0(dw\times[0,1))=\lambda$,
or equivalently, in terms of $\sigma$, we assume
\begin{equation}
\eqna{Burgersassumps}
\int_W \sigma(w,y)\,\lambda(dw)=1,\ y\in[0,1],
\end{equation}
and
\begin{equation}
\eqna{Burgersdensity1}
\int_0^1\sigma(w,z)\,dz=1,\ w\in W.
\end{equation}
We now state the main result we prove in this paper.
For notational convenience, in the following, and throughout the paper,
we use a notation such as $\mu(dw)=\nu(dw)$ to indicate
the equality of measures, $\mu(B)=\nu(B)$, for all $B\in{\cal B}(W)$.
\begin{thm}
\thma{Burgers}
There exists a unique pair of functions $y_C$ and $\mu_t(dw\times dz)$, 
where $y_C$ is a function of $(\gamma,t)\in\Delta_T$ taking values in
$[0,1]$, and $\mu_t(dw\times dz)$ is a function of $t\in [0,T]$
taking values in the probability measures on $W\times[0,1]$,
such that the following hold.
\begin{enumerate}
\item
$y_C((y_0,0),t)$ is non-decreasing in $y_0$,
$y_C((0,t_0),t)$ is non-increasing in $t_0$,
and $y_C(\gamma,t)$ is non-decreasing in $t$.
\item
$y_C(\gamma,t)$ and $\dsp \pderiv{y_C}{t}(\gamma,t)$ are continuous,
and for each $t\in[0,T]$, $y_C(\cdot,t):\ \Gamma_t\to[0,1]$ is surjective.
\item
For all bounded measurable $h:\ W\to\reals$,
$\dsp \int_W h(w) \mu_t(dw\times[y,1))$ is
Lipschitz continuous in $(y,t)\in [0,1]\times[0,T]$,
with Lipschitz constant uniform in $h$ satisfying
\begin{equation}
\eqna{Burgersprf3}
\sup_{w\in W} |h(w)|\le 1.
\end{equation}
More precisely,
\begin{equation}
\eqna{globalLipschitz}
\begin{array}{l}\dsp
\biggl|\int_W h(w)\mu_{t'}(dw\times[y',1)))
-\int_W h(w)\mu_{t}(dw\times[y,1))\biggr|
\le |y'-y|+M_We^{2C_WT}|t'-t|,
\end{array}
\end{equation}
for $h$ satisfying \eqnu{Burgersprf3}.
\item
The following equation of motion and initial and boundary conditions
hold.
\end{enumerate}
\begin{equation}
\eqna{U0}
y_C((y_0,t_0),t_0)=y_0\,,\ (y_0,t_0)\in\Gamma,\ \ \ 
\ \mbox{ and }\ 
\mu_0(dw\times dy)\ \mbox{ as in \eqnu{Burgersdefdensity},}
\end{equation}
\begin{equation}
\eqna{conservation}
\mu_t(dw\times[0,1))=\lambda(dw),\ \ t\in[0,T],
\end{equation}
\begin{equation}
\eqna{soliditycond}
\mu_t(W\times[y,1))=1-y,\ \ (y,t)\in[0,1]\times[0,T],
\end{equation}
\begin{equation}
\eqna{phic}
\begin{array}{l}\dsp
\mu_t(dw\times[y_C((y_0,t_0),t),1))
\\ \dsp
=\mu_{t_0}(dw\times[y_0,1))
-\int_{t_0}^t \int_{z\in[y_C((y_0,t_0),s),1)} w(z,s) 
\mu_s(dw\times dz)\,ds,
\\ \dsp
\ \ ((y_0,t_0),t)\in \Delta_T\,.
\end{array}
\end{equation}
\DDD
\end{thm}
Note that a substitution $y=y_C(\gamma,t)$ in \eqnu{soliditycond} implies
\begin{equation}
\eqna{varphiyC}
y_C(\gamma,t)=1-\mu_t(W\times [y_C(\gamma,t),1)),
\end{equation}
with which \eqnu{phic} and \eqnu{soliditycond} imply
\begin{equation}
\eqna{yc}
y_C(\gamma,t)=y_0+\int_{t_0}^t \int_{W\times[y_C(\gamma,s),1)}
 w(z,s) \mu_s(dw\times dz)\,ds.
\end{equation}

If $W$ is a countable set $W=\{w_1,w_2,\ldots\}$, denote
the distribution functions by
\begin{equation}
\eqna{mu2U}
 U_{\alpha}(y,t)=\mu_t(\{w_{\alpha}\}\times [y,1)). 
\end{equation}
Assume further that the functions 
$U_{\alpha}:\ [0,1]\times[0,T]\to\preals$
are in $C^1$.
Differentiating \eqnu{yc} by $t$ 
we reproduce \eqnu{fluidvelocity} in \secu{intro}.
Differentiating \eqnu{phic} by $t$, substituting \eqnu{fluidvelocity},
and then changing the notation from $y_C(\gamma,t)$ to $y$, 
we can eliminate the
dependence on initial/boundary parameter $\gamma$, and we reproduce
\eqnu{Burgers} in \secu{intro}.
With $\dsp \lambda(\{w_{\alpha}\})=r_{\alpha}$,
\eqnu{conservation} and \eqnu{soliditycond} respectively correspond to
\eqnu{Burgersconservation} and \eqnu{Burgerssoliditycond},
and the conditions \eqnu{rw} and \eqnu{Cw} imply \eqnu{Burgersrw}.
Thus \thmu{Burgers} contains a solution to the problem 
introduced in \secu{intro}.

In \thmu{Burgers} we claim differentiability for $y_C(\gamma,t)$ in $t$,
while we formulated \eqnu{phic} so that differentiability assumptions 
on $U_{\alpha}(y,t)$ or $\mu_t(dw\times [y,1))$ are absent.
In fact, at $(y,t)$ with $y=y_C((0,0),t)$, where the characteristic
curves starting at initial points $\gamma\in\Gamma_i$ 
and those starting at boundary points $\gamma\in\Gamma_b$ meet,
the differentiability with respect to variables which cross
the curve are lost in general. 
Loss of regularity across the characteristic curves is common 
for the quasilinear partial differential equations \cite{Bressan}.
In terms of \cite[\S 3.4]{Bressan},
we may therefore say that \thmu{Burgers} claims 
global existence of the Lipschitz solution
(broad solution which is Lipschitz continuous) to the system of
quasilinear partial differential equations \eqnu{Burgers},
where we extended the definition of Lipschitz solution
in \cite[\S 3.4]{Bressan},
to include the non-local (integration) terms,
and also generalized the notion of domain of determinancy
defined in \cite[\S 3.4]{Bressan}, which in the present case corresponds to
$\dsp \{(y,t)\in[0,1]\times\preals\mid y\ge y_C((0,0),t)\}$,
to the boundary condition dependent domain 
$\dsp \{(y,t)\in[0,1]\times\preals\mid y<y_C((0,0),t)\}$.
By formulating \thmu{Burgers} in terms of probability measures 
on $W\times[0,1]$ we also included 
uncountably many components parametrized by the evaporation rates $w$,
which are componentwise bounded but may be unbounded as a total fluid.

\section{Point process with last-arrival-time dependent intensity.}
\seca{pointprocess}

Let $N=N(t)$, $t\ge0$, be a non-decreasing, right-continuous, 
non-negative integer valued stochastic process
on a measurable space with $N(0)=0$,
and for each non-negative integer $k$ define its $k$-th arrival time $\tau_k$
by
\begin{equation}
\eqna{arrivaltime}
\tau_k=\inf\{t\ge0\mid N(t)\ge k\},\ \ k=1,2,\ldots,
\ \mbox{ and }\ \tau_0=0.
\end{equation}
The arrival times $\tau_k$ are non-decreasing in $k$, 
because $N$ is non-decreasing,
and since $N$ is also right-continuous, 
the arrival times are stopping times;
if we denote the associated filtration by 
$\dsp {\cal F}_t=\sigma[N(s),\ s\le t]$,
then $\{\tau_k\le t\} \in {\cal F}_t$, $t\ge 0$.

Let $\omega$ be a non-negative valued bounded continuous function 
of $(s,t)$ for $0\le s\le t$, and for $k=1,2,\ldots$ assume that
\begin{equation}
\eqna{intensity}
\begin{array}{l}\dsp
\prb{t< \tau_k\mid {\cal F}_{\tau_{k-1}}}
=\exp(-\int_{\tau_{k-1}}^t \omega(\tau_{k-1},u)\,du)
\ \mbox{ on }\ t\ge \tau_{k-1}\,.
\end{array}
\end{equation}
In particular, \eqnu{intensity} with $k=1$ implies
\begin{equation}
\eqna{intensity0}
\prb{N(t)=0}=\prb{\tau_1>t}=\exp(-\int_0^t \omega(0,u)\,du),\ \ t\ge0.
\end{equation}

Note that the function $\omega$ has different dependence on the variables
from the evaporation rate function $w$ in the other sections of this paper.
(We will relate $\omega$ to $w$ by \eqnu{w2omega} in \secu{flows},
namely, we will introduce an intensity function as a composite function
of the evaporation rate function and a flow.)
If $\omega$ is independent of the first variable, then
\eqnu{intensity} implies that $N$ is the (inhomogeneous) Poisson process
with intensity function $\omega$.
We are considering a generalization of the Poisson process
such that the intensity function depends on the latest arrival time.

The remainder of this section is devoted to 
basic formulas to be used in this paper.
For a continuously differentiable function $f$ vanishing at $\infty$,
\eqnu{intensity}, with integration by parts and the Fubini's theorem,
implies
\begin{equation}
\eqna{stoppingtimeFubini}
\begin{array}{l}\dsp
\int_{\tau_{k-1}}^{\infty} f(t)\,\omega(\tau_{k-1},t)\,
\exp(-\int_{\tau_{k-1}}^t \omega(\tau_{k-1},s)\,ds)\,dt
\\ \dsp \phantom{}
=\int_{\tau_{k-1}}^{\infty} f'(t)\,
\exp(-\int_{\tau_{k-1}}^t \omega(\tau_{k-1},s)\,ds)\,dt+f(\tau_{k-1})
\\ \dsp \phantom{}
=\int_{\tau_{k-1}}^{\infty} f'(t)\,
\prb{t< \tau_k\mid {\cal F}_{\tau_{k-1}}}\,dt+f(\tau_{k-1})
\\ \dsp \phantom{}
=\EE{\int_{\tau_{k-1}}^{\infty}
 f'(t)\chrfcn{t< \tau_k}\,dt\mid
{\cal F}_{\tau_{k-1}}}+f(\tau_{k-1})
\\ \dsp \phantom{}
=\EE{\int_{\tau_{k-1}}^{\tau_k} f'(t)\,dt\mid {\cal F}_{\tau_{k-1}}}
+f(\tau_{k-1})
\\ \dsp \phantom{}
=\EE{f(\tau_k)\mid {\cal F}_{\tau_{k-1}}},\ \ k=1,2,\ldots.
\end{array}
\end{equation}
Approximating by a series of smooth functions, 
\eqnu{stoppingtimeFubini} holds for any $f\in L_0(\preals)$,
where $L_0(\preals)$ is 
the space of bounded measurable functions 
$f:\ \preals\to\reals$ vanishing at infinity,
equipped with the supremum norm.

For $t\ge t_0$ put
\begin{equation}
\eqna{numberdepPoissonprf11}
\Omega(t_0,t)=\int_{t_0}^t \omega(t_0,u)\,du,
\end{equation}
and define a linear map $A_{\omega}:\ L_0(\preals)\to L_0(\preals)$ by
\begin{equation}
\eqna{stoppingtimeopA}
(A_{\omega}f)(t)=\int_t^{\infty} f(u)\, \omega(t,u)\, e^{-\Omega(t,u)}\,du.
\end{equation}
Then \eqnu{stoppingtimeFubini} implies
\begin{equation}
\eqna{stoppingtimeop}
\EE{f(\tau_k)\mid {\cal F}_{\tau_{k-1}}}=(A_{\omega}f)(\tau_{k-1}),
\ \ f\in L_0(\preals).
\end{equation}
By induction and $\tau_0=0$ we have
\begin{equation}
\eqna{stoppingtimeexpectation}
\begin{array}{l}\dsp
\EE{f(\tau_k)}
=\EE{ \EE{ \cdots \EE{ \EE{f(\tau_k)\mid {\cal F}_{\tau_{k-1}}}
 \mid {\cal F}_{\tau_{k-2}} } \cdots \mid {\cal F}_{\tau_1 }} }
\\ \dsp \phantom{\EE{f(\tau_k)}}
=(A_{\omega}^k\,f)(0)
\\ \dsp \phantom{\EE{f(\tau_k)}}
=\int_{0\le u_1\le u_2\le \cdots\le u_k<\infty} f(u_k)\,
\prod_{i=1}^k \omega(u_{i-1},u_i)\,e^{-\Omega(u_{i-1},u_i)}
\, du_i\,,
\end{array}
\end{equation}
where we put $u_0=0$ to simplify notations.

For example, by choosing $\dsp f(u)=\chrfcn{u\le t}$, 
\eqnu{stoppingtimeop} implies
\begin{equation}
\eqna{Poissonprob1}
\prb{\tau_k\le t\mid {\cal F}_{\tau_{k-1}}}
=(1-e^{-\Omega(\tau_{k-1},t)})\,\chrfcn{\tau_{k-1}\le t}\,.
\end{equation}
Then 
\eqnu{Poissonprob1} and
\eqnu{stoppingtimeexpectation} with
$\dsp f(u)=e^{-\Omega(u,t)}\,\chrfcn{u\le t}$
imply
\begin{equation}
\eqna{Poissonprb}
\begin{array}{l}\dsp
\prb{N(t)=k}=\prb{\tau_k\le t<\tau_{k+1}}
=\EE{\chrfcn{\tau_k\le t}\,(1-\prb{\tau_{k+1}\le t\mid {\cal F}_{\tau_k}})\,}
\\ \dsp\phantom{\prb{N(t)=k}}
=\EE{f(\tau_k)}
=(A_{\omega}^k\,f)(0),\ \ k\in\pintegers,\ t>0.
\end{array}
\end{equation}
Hence, as in the last line of \eqnu{stoppingtimeexpectation},
\begin{equation}
\eqna{Poissonprbexplicit}
\begin{array}{l}\dsp
\prb{N(t)=k}=(A_{\omega}^k f)(0)
\\ \dsp\phantom{}
=\left\{\begin{array}{ll}\dsp
\int_{0\le u_1\le u_2\le \cdots\le u_k\le t} e^{-\Omega(u_k,t)}\,
\prod_{i=1}^k  \omega(u_{i-1},u_i)\,e^{-\Omega(u_{i-1},u_i)}
\, du_i\,, & k=1,2,\ldots,
\\ \dsp e^{-\Omega(0,t)}, & k=0. \end{array} \right.
\end{array}
\end{equation}
In particular, $\prb{N(t)\ge 0}=1$ implies a sum rule
\begin{equation}
\eqna{Poissonsumrule}
\begin{array}{l}\dsp
e^{-\Omega(0,t)}+
\sum_{k=1}^{\infty}
\int_{0\le u_1\le u_2\le \cdots\le u_k\le t} e^{-\Omega(u_k,t)}\,
\prod_{i=1}^k  \omega(u_{i-1},u_i)\,e^{-\Omega(u_{i-1},u_i)}
\, du_i=1,
\\ \dsp
t>0,
\end{array}
\end{equation}
where $u_0=0$, as in \eqnu{stoppingtimeexpectation}.

Similarly, given $s$ and $t$ satisfying $0\le s<t$,
the probability that there is no arrival in the interval $(s,t]$ is
\begin{equation}
\eqna{Poissondecay}
\begin{array}{l}\dsp
\prb{N(t)=N(s)}
=\sum_{k=0}^{\infty} \prb{N(t)=N(s)=k}
=\sum_{k=0}^{\infty} \prb{\tau_k\le s,\ t<\tau_{k+1}}
\\ \dsp \phantom{\prb{N(t)=N(s)}}
=\sum_{k=0}^{\infty} \EE{%
\chrfcn{\tau_k\le s}\,(1-\prb{\tau_{k+1}\le t\mid{\cal F}_{\tau_k}})\,}
\\ \dsp \phantom{\prb{N(t)=N(s)}}
=\sum_{k=0}^{\infty} \EE{\chrfcn{\tau_k\le s}\,e^{-\Omega(\tau_k,t) }}.
\end{array}
\end{equation}
With $\dsp f(u)=\chrfcn{u\le s}\,e^{-\Omega(u,t)}$
in \eqnu{stoppingtimeexpectation}, we also have an explicit formula
\begin{equation}
\eqna{Poissondecayexplicit}
\begin{array}{l}\dsp
\prb{N(t)=N(s)=k}= \EE{\chrfcn{\tau_k\le s}\,e^{-\Omega(\tau_k,t) }}
=(A_{\omega}^k\,f)(0)
\\ \dsp 
=\left\{\begin{array}{ll}\dsp e^{-\Omega(0,t)}, & k=0,\\ \dsp
\int_{0\le u_1\le u_2\le \cdots\le u_k\le s} e^{-\Omega(u_k,t)}\,
\prod_{i=1}^k \omega(u_{i-1},u_i)\,e^{-\Omega(u_{i-1},u_i)}\,du_i
\,, & k=1,2,\ldots, \end{array} \right.
\end{array}
\end{equation}
for $t\ge s>0$, where $u_0=0$, as in \eqnu{stoppingtimeexpectation}.
The following property relates the $s$ and $t$ dependencies
of the quantity in \eqnu{Poissondecay}. Note that the explicit formula
\eqnu{Poissondecay} implies that this quantity is $C^1$ in $s$ and $t$.
\begin{prp}
\prpa{st_dep}
For $k=1,2,\ldots$,
\begin{equation}
\eqna{st_dep}
\begin{array}{l}\dsp
\pderiv{}{t}\prb{N(t)=N(s)=k}
=-\int_0^s \omega(u,t)\, \pderiv{}{u}\prb{N(t)=N(u)=k}\,du,
\\ \dsp
0\le s<t.
\end{array}
\end{equation}
\DDD
\end{prp}
\begin{prf}
First we prove
\begin{equation}
\eqna{st_dep_prf1}
\EE{f(\tau_k)\,g(\tau_k)\,\chrfcn{\tau_k\le s}}
=\int_0^s f(u)\,Q'(u)\,du
\end{equation}
for locally bounded and measurable $f$ and $g$ such that
\begin{equation}
\eqna{st_dep_prf2}
Q(s):=\EE{g(\tau_k)\,\chrfcn{\tau_k\le s}}
\end{equation}
is absolutely continuous with respect to the Lebesgue measure
(so that the derivative $Q'$ almost surely exists).
Approximating by a series of smooth functions, 
it suffices to prove \eqnu{st_dep_prf1} for $f\in C^1$.
By Fubini's theorem and partial integration,
and noting that $\tau_k>0$ for $k>0$ implies $Q(0)=0$,
\[ \begin{array}{l}\dsp
\EE{f(\tau_k)\,g(\tau_k)\,\chrfcn{\tau_k\le s}}
=\EE{\biggl(-\int_{\tau_k}^s f'(u)\,du+f(s)\biggr)\,g(\tau_k)
\,\chrfcn{\tau_k\le s}}
\\ \dsp \phantom{\EE{f(\tau_k)\,g(\tau_k)\,\chrfcn{\tau_k\le s}}}
=f(s)\,Q(s)-\int_0^sf'(u)\,\EE{\chrfcn{\tau_k\le u}\,g(\tau_k)}\,du
\\ \dsp \phantom{\EE{f(\tau_k)\,g(\tau_k)\,\chrfcn{\tau_k\le s}}}
=f(s)\,Q(s)-\int_0^s f'(u)\,Q(u)\,du
\\ \dsp \phantom{\EE{f(\tau_k)\,g(\tau_k)\,\chrfcn{\tau_k\le s}}}
=\int_0^s f(u)\,Q'(u)\,du.
\end{array} \]
Thus \eqnu{st_dep_prf1} is proved.

Now for a positive integer $k$,
let $f(u)=\omega(u,t)$ and $\dsp g(u)=e^{-\Omega(u,t)}$ in
\eqnu{st_dep_prf1}. 
Note that for this choice \eqnu{Poissondecayexplicit} implies
\[ 
Q(s)=\EE{e^{-\Omega(\tau_k,t)}\,\chrfcn{\tau_k\le s}}
=\prb{N(t)=N(s)=k}.
\]
Then \eqnu{st_dep_prf1} implies
\[ \begin{array}{l} \dsp
\pderiv{}{t}\prb{N(t)=N(s)=k}
=-\EE{\omega(\tau_k,t)\,e^{-\Omega(\tau_k,t)}\,\chrfcn{\tau_k\le s} }
\\ \dsp \phantom{\pderiv{}{t}\prb{N(t)=N(s)=k}}
=-\EE{f(\tau_k)\,g(\tau_k)\,\chrfcn{\tau_k\le s} }
=-\int_0^s f(u)\,Q'(u)\,du
\\ \dsp \phantom{\pderiv{}{t}\prb{N(t)=N(s)=k}}
=-\int_0^s \omega(u,t)\,\pderiv{}{u} \EE{N(t)=N(u)=k}\,du,
\end{array}\]
which proves \eqnu{st_dep}.
\QED
\end{prf}

Note that, for example, the explicit formula in
\eqnu{Poissondecayexplicit} depends on the intensity function $\omega$
at times before $s$, which implies that the process $N$ is not of
independent increment, hence, in particular, is not a Poisson process.
If, on the other hand, $\omega$ is independent of its first variable,
put $\dsp\tilde{\omega}(t)=\omega(s,t)$. Then
\eqnu{numberdepPoissonprf11} implies 
$\dsp \Omega(t_0,t)=\int_{t_0}^t \tilde{\omega}(s)\,ds$, and
\eqnu{Poissondecayexplicit} is simplified as
\begin{equation}
\eqna{Poissondecayindepcase}
\begin{array}{l}\dsp
\prb{N(t)=N(s)}=\biggl(1+\sum_{k\ge1}
\int_{0\le u_1\le u_2\le \cdots\le u_k\le s}
\prod_{i=1}^k \tilde{\omega}(u_i)\,
\prod_{i=1}^k du_i\biggr)\,e^{-\Omega(0,t)}
\\ \dsp \phantom{\prb{N(t)=N(s)}}
=\biggl(1+\sum_{k\ge1}\frac1{k!}\,\Omega(0,s)^k\biggr)\,e^{-\Omega(0,t)}
\\ \dsp \phantom{\prb{N(t)=N(s)}}
=e^{\Omega(0,s)}\times e^{-\Omega(0,t)}
=e^{-\Omega(s,t)},
\end{array}
\end{equation}
where we used an elementary formula proved by induction in $k$,
\begin{equation}
\eqna{integralopcompactness}
\int_{0\le u_1\le u_2\le \cdots\le u_k\le s}
\prod_{i=1}^k f(u_i)du_1\,du_2\ldots du_k
=\frac1{k!}\biggl(\int_0^s f(v)dv\biggr)^k,
\ \ s\ge0,\ k=1,2,\ldots,
\end{equation}
valid for any integrable function $f:\ \reals\to\reals$.
The simple result \eqnu{Poissondecayindepcase}
reproduces a formula for the (inhomogeneous)
Poisson process with independent increments.
In the general case of processes we consider in this paper,
such simple relations to Poisson processes or Poisson distributions
are lost and the properties of the processes become complicated.

\section{Flows and construction of solution.}
\seca{flows}

The key quantities for the solution to the functional equations in 
\thmu{Burgers} are the characteristic curves $y_C$ and the associated 
measure $\dsp \varphi(dw,\gamma,t)=\mu_t(dw\times[y_C(\gamma,t),1))$.
We will find $y_C$ as a unique solution to a non-linear map
on a space $\Theta_T$ of flows, a non-decreasing function in time $t$ and in
initial/boundary points $\gamma\in \Gamma$.
To simplify the definition of $\Theta_T$ we first define a total order
$\succeq$ on $\Gamma$ by
\begin{equation}
\eqna{Gammaorder}
s\le t,\ z\le y\ \Leftrightarrow\ (0,T)\succeq (0,t)\succeq (0,s)\succeq (0,0)
\succeq (z,0)\succeq (y,0)\succeq (1,0).
\end{equation}
We now define the set of flows $\Theta_T$ on $[0,1]\times[0,T]$ by
\begin{equation}
\eqna{class_of_conti_flows}
\begin{array}{l}\dsp
\Theta_T:=\{\theta:\ \Delta_T\to[0,1]\mid 
\theta((y_0,t_0),t_0)=y_0,\ (y_0,t_0)\in\Gamma_T,\ \mbox{ continuous,
 }\\ \dsp \phantom{\Theta_T:=\{} \mbox{ %
surjective and non-increasing in $\gamma$ for each $t$,
 }\\ \dsp \phantom{\Theta_T:=\{} \mbox{ %
non-decreasing in $t$ for each $\gamma$
 }\}.
\end{array}
\end{equation}
For example,
\begin{equation}
\eqna{Gamma_infinity}
\theta((1,0),t)=1,\ \ t\in[0,T],\ \theta\in\Theta_T.
\end{equation}

Let $W$, the set of evaporation rates, be as in \thmu{Burgers},
and let $\theta\in\Theta_T$.
For each $w\in W$ and $z\in[0,1)$ define 
$\omega=\omega_{\theta,w,z}$, a non-negative valued continuous function 
of $(s,t)$ satisfying $0\le s\le t\le T$, by
\begin{equation}
\eqna{w2omega}
\omega_{\theta,w,z}(s,t)=\left\{\begin{array}{ll}\dsp w(\theta((z,0),t),t), 
& \mbox{if }\ s=0,\\ \dsp w(\theta((0,s),t),t),& \mbox{if }\ s>0.
\end{array}\right.
\end{equation}
Note that $\omega_{\theta,w,z}$ is independent of $z$ if $s>0$.
Let $\dsp \{N_{\theta,w,z}\mid z\in[0,1),\ w\in W\}$ be a set of 
processes, with each $N_{\theta,w,z}$ being a point process
$N$ introduced in \secu{pointprocess} with the intensity function
in \eqnu{intensity} determined by $\omega=\omega_{\theta,w,z}$.
The quantity in \eqnu{numberdepPoissonprf11} for
the choice \eqnu{w2omega} is
\begin{equation}
\eqna{omegaOmega}
\begin{array}{l}\dsp
\Omega_{\theta,w,z}(0,t)=\int_0^t w(\theta((z,0),u),u)\,du,
\\ \dsp
\Omega_{\theta,w}(s,t)=\int_s^t w(\theta((0,s),u),u)\,du,
\ \ 0<s\le t.
\end{array}
\end{equation}

Let $\mu_0$ be as in \thmu{Burgers}, and
define a function $\varphi_{\theta}(dw,\gamma,t)$ 
on $(\gamma,t)\in \Delta_T$ taking values in the measures on $W$,
by
\begin{equation}
\eqna{varphi_theta}
\begin{array}{l}\dsp
\varphi_{\theta}(dw,\gamma,t)=\int_{z\in[y_0,1)}
 \prb{N_{\theta,w,z}(t)=N_{\theta,w,z}(t_0)}\,\mu_0(dw\times dz),
\\ \dsp
\ \gamma=(y_0,t_0)\in\Gamma,\ (\gamma,t)\in\Delta_T\,.
\end{array}
\end{equation}
The explicit form for \eqnu{varphi_theta} is simple for
$\dsp\gamma=(y_0,0)\in \Gamma_i$, because $\dsp N_{\theta,w,z}(0)=0$,
and \eqnu{Poissondecayexplicit} with $k=0$ imply
\begin{equation}
\eqna{varpphi_i}
\varphi_{\theta}(dw,(y_0,0),t)
=\int_{z\in[y_0,1)} e^{-\Omega_{\theta,w,z}(0,t)}\,\mu_0(dw\times dz).
\end{equation}
For $\varphi_{\theta}$ in \eqnu{varphi_theta}
define $\dsp \pderiv{\varphi_{\theta}}{\gamma}$,
a measure valued function on $\Delta_T$, by
\begin{equation}
\eqna{dphidgamma}
\pderiv{\varphi_{\theta}}{\gamma}(dw,\gamma,t)
=\left\{\begin{array}{ll} \dsp -\pderiv{\varphi_{\theta}}{z}(dw,(z,0),t),
& \mbox{ if }\ \gamma=(z,0)\in\Gamma_i\,, \\ \dsp
\pderiv{\varphi_{\theta}}{u}(dw,(0,u),t),
& \mbox{ if }\ \gamma=(0,u)\in\Gamma_b\,. \end{array} \right.
\end{equation}
We keep non-negativity of the defined measure in 
determining the sign.
Explicit calculation of the derivative at $\gamma\in\Gamma_i$ is 
straightforward from \eqnu{varpphi_i} and \eqnu{Burgersdefdensity}.
The derivative at $\gamma=(0,u)\in\Gamma_b$ is also 
calculated explicitly using \eqnu{Poissondecayexplicit} and 
\eqnu{varphi_theta}, which is
\begin{equation}
\eqna{dphidu}
\begin{array}{l}\dsp
\pderiv{\varphi_{\theta}}{u}(dw,(0,u),t)
\\ \dsp \phantom{}
=\int_{z\in[0,1)}\biggl( w(\theta((z,0),u)\,e^{-\Omega_{\theta,w,z}(0,u)}
\\ \dsp \phantom{=\int}
+\sum_{k\ge2} \int_{0< u_1\le u_2\le \cdots\le u_{k-1}\le u} 
w(\theta((z,0),u_1)\,e^{-\Omega_{\theta,w,z}(0,u_1)}\,du_1
\\ \dsp \phantom{\int=+\sum_{k\ge2} \int}
\prod_{i=2}^{k-1} 
 \biggl( w(\theta((0,u_{i-1}),u_i)
\,e^{-\Omega_{\theta,w}(u_{i-1},u_i)}\,du_i\biggr)
\\ \dsp \phantom{\int=+\sum_{k\ge2} \int}
w(\theta((0,u_{k-1}),u)\,e^{-\Omega_{\theta,w}(u_{k-1},u)}
\biggr)\,e^{-\Omega_{\theta,w}(u,t)}\,\mu_0(dw\times dz),
\end{array}
\end{equation}
where we also used the notations \eqnu{w2omega} and \eqnu{omegaOmega}
to make the $z$ and $u$ dependence explicit.
Note that $w\in W$ are non-negative functions,
hence $\dsp e^{-\Omega_{\theta,w}(u,t)}\le 1$. This and the sum rule
\eqnu{Poissonsumrule} (with the replacements $t=u$ and $k=k'-1$),
with \eqnu{rw}, \eqnu{Burgersdefdensity}, and \eqnu{Burgersdensity1} imply
\[
\pderiv{\varphi_{\theta}}{u}(W,(0,u),t)
\le \int_{B\times[0,1)} \wnorm{w}\,\mu_0(dw\times dz)
= \int_B \wnorm{w} \lambda(dw)=M_W<\infty,
\]
hence, $\dsp \pderiv{\varphi_{\theta}}{u}$ is well-defined.

For $f:\ \Gamma_T\to\reals$ and a Borel subset $A\subset \Gamma_T$\,,
define $\dsp\int_A f(\gamma)\,d\gamma$, a line integral on $\Gamma_T$,
by
\begin{equation}
\eqna{Gammaintegral}
\int_A f(\gamma)\,d\gamma=\int_{A_i} f(z,0)\,dz+\int_{A_b} f(0,u)\,du,
\end{equation}
where,
$\dsp A_i=\{z\in[0,1)\mid (z,0)\in A\}$ and
$\dsp A_b=\{u\in[0,T]\mid (0,u)\in A\}$,
and the integration in the right hand side of \eqnu{Gammaintegral} are 
the standard one dimensional integrations.
\begin{prp}
\prpa{finiteevaporation}
It holds that
\begin{equation}
\eqna{finiteevaporation}
0\le \int_{\gamma'\in\Gamma_t} \int_W 
 \wnorm{w}\,\pderiv{\varphi_{\theta}}{\gamma}(dw,\gamma',t)\,d\gamma'
 \le M_W\,e^{2C_Wt},
\end{equation}
for all $t\in[0,T]$.
\DDD
\end{prp}
\begin{rem}
Note the extra $\wnorm{w}$ in the 
integrand of \eqnu{finiteevaporation}.
We are allowing $W$ to contain unbounded functions $w$,
so the finiteness of \eqnu{finiteevaporation} is harder than
that of \eqnu{dphidu}.
We use the condition \eqnu{Cw} as well as \eqnu{rw}
to prove \eqnu{finiteevaporation}.

\prpu{finiteevaporation} implies, in particular, that the integration
\eqnu{currentconservation} introduced later, is uniformly bounded 
in $\dsp (\gamma,t)\in\Delta_T$ and $\dsp B\in{\cal B}(W)$.
\DDD
\end{rem}
\begin{prfof}{\protect\prpu{finiteevaporation}}
By the definitions
\eqnu{dphidgamma}
and
\eqnu{Gammaintegral},
\begin{equation}
\eqna{finiteevaporationprf1}
\begin{array}{l}\dsp
 \int_{\gamma'\in\Gamma_t} \int_W 
 \wnorm{w}\,\pderiv{\varphi_{\theta}}{\gamma}(dw,\gamma',t)\,d\gamma'
\\ \dsp {}
=
-\int_0^1\int_W \wnorm{w} \pderiv{\varphi_{\theta}}{z}(dw,(z,0),t)\,dz
+\int_0^t\int_W \wnorm{w} \pderiv{\varphi_{\theta}}{u}(dw,(0,u),t)\,du.
\end{array}
\end{equation}
The first term on the right hand side is explicitly calculated using 
\eqnu{varpphi_i} and \eqnu{Burgersdefdensity}.
Using also $\dsp e^{-\Omega_{\theta,w,z}(0,t)}\le 1$, 
\eqnu{Burgersdensity1}, and \eqnu{rw}, we have an estimate
\begin{equation}
\eqna{finiteevaporationprf2}
\begin{array}{l}\dsp
0\le
-\int_0^1\int_W \wnorm{w} \pderiv{\varphi_{\theta}}{z}(dw,(z,0),t)\,dz
=\int_{W\times [0,1)} \wnorm{w} e^{-\Omega_{\theta,w,z}(0,t)}
\,\mu_0(dw\times dz)
\\ \dsp \phantom{0}
\le M_W<\infty.
\end{array}
\end{equation}
Next, for $t\ge s\ge 0$ and $w\in W$, put
\begin{equation}
\eqna{wstar}
\tilde{\Omega}_w(s,t)=\int_s^t w(1,u)\,du.
\end{equation}
Then the condition \eqnu{Cw} and
the fact that $\theta$ takes values in $[0,1]$ imply
\begin{equation}
\eqna{finiteevaporationprf3}
\begin{array}{l}\dsp
|w(\theta(\gamma,t),t)-w(1,t)|\le C_W\,,
\\ \dsp
|\Omega_{\theta,w}(s,t)-\tilde{\Omega}_w(s,t)|\le C_W\,(t-s),
\ \mbox{ and }\ 
|\Omega_{\theta,w,z}(0,t)-\tilde{\Omega}_w(0,t)|\le C_W\,t,
\\
(\gamma,t)\in\Delta_T\,,\ 0< s\le t,\  w\in W.
\end{array}
\end{equation}
Using \eqnu{finiteevaporationprf3} and \eqnu{wstar} in
\eqnu{dphidu}, and then using
\eqnu{Burgersdefdensity} and \eqnu{Burgersdensity1},
we have an estimate
\[ \begin{array}{l}\dsp
\int_0^t\int_W \wnorm{w} \pderiv{\varphi_{\theta}}{u}(dw,(0,u),t)\,du
\\ \dsp{}
\le\int_W \wnorm{w}\,e^{-\tilde{\Omega}_w(0,t)+C_W\,t}\,
\int_0^t(w(1,u)+C_W)\,
\biggl(1
\\ \dsp \phantom{\le\int_W\wnorm{w}\,}
+\sum_{k\ge 1} \int_{0<u_1\le \cdots\le u_{k}\le u}
\prod_{i=1}^k ((w(1,u_i)+C_W)\,du_i)\biggr)
\,du\,\lambda(dw).
\end{array}\]
The estimate is now reduced to that for the Poisson processes,
and \eqnu{integralopcompactness} implies
\[ \begin{array}{l}\dsp
\int_0^t\int_W \wnorm{w} \pderiv{\varphi_{\theta}}{u}(dw,(0,u),t)\,du
\\ \dsp{}
\le
\int_W \wnorm{w}e^{-\tilde{\Omega}_w(0,t)+C_Wt}
\int_0^t e^{\tilde{\Omega}_w(0,u)+C_Wu}(w(1,u)+C_Wu)\,du\,\lambda(dw)
\\ \dsp{}
=
\int_W \wnorm{w}e^{-\tilde{\Omega}_w(0,t)+C_Wt}
\biggl[e^{\tilde{\Omega}_w(0,u)+C_Wu}\biggr]_{u=0}^{u=t}\lambda(dw)
\\ \dsp{}
\le e^{2C_Wt}\int_W\wnorm{w}\lambda(dw) \le M_W\,e^{2C_Wt}.
\end{array}\]
This proves \eqnu{finiteevaporation}.
\QED
\end{prfof}
\begin{prp}
\prpa{currentconservation}
It holds that
\begin{equation}
\eqna{currentconservation}
\begin{array}{l}\dsp
\pderiv{\varphi_{\theta}}{t}(B,\gamma,t)
=-\int_{\gamma\succeq\gamma'} \int_B 
 w(\theta(\gamma',t),t)\,\pderiv{\varphi_{\theta}}{\gamma}(dw,\gamma',t)
\,d\gamma',
\\ \dsp
(\gamma,t)\in\Delta_T,\ B\in{\cal B}(W),
\end{array}
\end{equation}
where $\gamma\succeq\gamma'$ is defined in \eqnu{Gammaorder}.
\DDD
\end{prp}
\begin{rem}
If $W$ consists of functions with no spatial dependence, 
namely, if $w(y,t)=w(1,t)$, then the factor $w$ in
the integrand of the right hand side
of \eqnu{currentconservation} is constant for $\gamma'$ integration,
and we have integration after differentiation,
so that the right hand side is simplified as 
$\dsp -w(1,t) \varphi_{\theta}(B,\gamma,t)$,
and the equation is solved easily,
as remarked below \eqnu{Burgersspaceindepfluidvelocity} in \secu{intro}.
\DDD\end{rem}
\begin{prfof}{\protect\prpu{currentconservation}}
Consider first the case $\gamma=(y_0,0)\in\Gamma_i$\,.
The explicit form \eqnu{varpphi_i},
together with the definitions
\eqnu{dphidgamma},
\eqnu{Gammaintegral}, and
\eqnu{Burgersdefdensity},
implies
\[ \begin{array}{l}\dsp
-\int_{\gamma\succeq\gamma'} \int_B 
 w(\theta(\gamma',t),t)\,\pderiv{\varphi_{\theta}}{\gamma}(dw,\gamma',t)
\,d\gamma'
\\ \dsp {}
=
-\int_{z\in[y_0,1)} \int_B w(\theta((z,0),t),t)\,
e^{-\Omega_{\theta,w,z}(0,t)}\,\sigma(w,z)\lambda(dw)\,dz
\\ \dsp {}
=
-\int_{z\in[y_0,1)}\int_B w(\theta((z,0),t),t)\,
e^{-\Omega_{\theta,w,z}(0,t)}\,\mu_0(dw\times dz)
=
\pderiv{\varphi_{\theta}}{t}(B,(y_0,0),t),
\end{array}
\]
which proves \eqnu{currentconservation} for $\gamma\in\Gamma_i$\,.

To prove \eqnu{currentconservation} for $\gamma\in \Gamma_b$\,,
put, for $k\in\pintegers$,
\begin{equation}
\eqna{currentconservationprf1}
\varphi\ss{k}_{\theta}(dw,(y_0,t_0),t)=\int_{z\in[y_0,1)}
 \prb{N_{\theta,w,z}(t)=N_{\theta,w,z}(t_0)=k}\,\mu_0(dw\times dz).
\end{equation}
Then \eqnu{varphi_theta} implies 
\begin{equation}
\eqna{currentconservationprf2}
\varphi_{\theta}=\sum_{k=0}^{\infty} \varphi\ss{k}_{\theta}.
\end{equation}
We will prove, for $\gamma=(0,t_0)\in \Gamma_b$
\begin{equation}
\eqna{currentconservationprf3}
\begin{array}{l}\dsp
\pderiv{\varphi\ss{k}_{\theta}}{t}(B,\gamma,t)
=-\int_{\gamma\succeq\gamma'} \int_B  w(\theta(\gamma',t),t)
\,\pderiv{\varphi\ss{k}_{\theta}}{\gamma}(dw,\gamma',t)\,d\gamma',
\\ \dsp
(\gamma,t)\in\Delta_T,\ B\in{\cal B}(W),
\end{array}
\end{equation}
for all $k\in\pintegers$, where differentiation and integration
with respect to $\gamma$ are defined in accordance with
\eqnu{dphidgamma} and \eqnu{Gammaintegral}.
Then \eqnu{currentconservationprf2} and 
\eqnu{currentconservationprf3} prove
\eqnu{currentconservation}.
The changes in the order of series and integration and differentiation
causes no problem, because all the terms and integrands are non-negative
and the results of summation and integration are bounded by
\prpu{finiteevaporation}.

Consider first the case $k=0$.
Then \eqnu{Poissondecayexplicit} implies
\begin{equation}
\eqna{varphi0}
\varphi\ss{0}_{\theta}(dw,(y_0,t_0),t)
=\int_{z\in[y_0,1)}e^{-\Omega_{\theta,w,z}(0,t)}\,\mu_0(dw\times dz).
\end{equation}
Note that this is independent of $t_0$.
Hence \eqnu{dphidgamma} and \eqnu{Burgersdefdensity} imply
\begin{equation}
\eqna{currentconservationprf4}
\pderiv{\varphi\ss0_{\theta}}{\gamma}(dw,\gamma',t)=
\left\{\begin{array}{ll} 
e^{-\Omega_{\theta,w,z}(0,t)}\,\sigma(w,z)\lambda(dw),& 
\mbox{ if } \gamma'=(z,0)\in\Gamma_i\,, \\ \dsp
0, & \mbox{ if } \gamma'=(0,u)\in\Gamma_b\,. \end{array}\right.
\end{equation}

For the case $\dsp\gamma=(0,t_0)\in\Gamma_b$\,,
\[
\gamma\succeq\gamma'\ \Leftrightarrow
\ \gamma'=(0,u),\ 0\le u\le t_0,
\ \mbox{ or }\ \gamma'=(z,0),\ 0\le z\le 1.
\]
The contribution, however, to \eqnu{Gammaintegral} from the integration
along $\Gamma_b$ vanishes because the integrand 
\eqnu{currentconservationprf4} is $0$ on $\Gamma_b$. Hence, 
\eqnu{currentconservationprf4},
\eqnu{Gammaintegral}, \eqnu{Burgersdefdensity}, and \eqnu{varphi0} imply
\[\begin{array}{l}\dsp
\int_{\gamma\succeq\gamma'} \int_B  w(\theta(\gamma',t),t)
\,\pderiv{\varphi\ss{0}_{\theta}}{\gamma}(dw,\gamma',t)
\,d\gamma'
\\ \dsp
=\int_{z\in[0,1)} \int_B  w(\theta((z,0),t),t)\,
 e^{-\Omega_{\theta,w,z}(0,t)}\,\sigma(w,z)\lambda(dw)\,dz
\\ \dsp
=\int_{B\times[0,1)}
 w(\theta((z,0),t),t)\,e^{-\Omega_{\theta,w,z}(0,t)}\,\mu_0(dw\times dz)
\\ \dsp
=-\pderiv{\varphi\ss{0}_{\theta}}{t}(B,\gamma,t),
\end{array}\]
which proves \eqnu{currentconservationprf3} for 
$k=0$ and $\gamma\in\Gamma_b$\,.

To consider the case $k>0$ and $\gamma=(0,t_0)\in\Gamma_b$,
\eqnu{currentconservationprf1}, \prpu{st_dep}, and
\eqnu{w2omega} imply
\[ \begin{array}{l}\dsp
\pderiv{\varphi\ss{k}_{\theta}}{t}(B,(0,t_0),t)
\\ \dsp
=-\int_{B\times [0,1)}
\int_0^{t_0} w(\theta(0,u),t),t)
\, \pderiv{}{u}\prb{N_{\theta,w,z}(t)=N_{\theta,w,z}(u)=k}\,du
\,\mu_0(dw\times dz)
\\ \dsp
=-\int_0^{t_0}\int_{w\in B} 
 w(\theta(0,u),t),t)\, \pderiv{\varphi\ss{k}_{\theta}}{u}(dw,(0,u),t)\,du,
\ \ 0< t_0<t.
\end{array}\]
The contribution to \eqnu{Gammaintegral} from $\gamma'\in\Gamma_i$
vanishes for the case $k>0$ because
$N_{\theta,w,z}(0)=0$ and \eqnu{currentconservationprf1}
then imply
\begin{equation}
\eqna{currentconservationprf5}
\varphi\ss{k}_{\theta}(dw,\gamma,t)=0,
\ \ \gamma\in\Gamma_i\,,\ t\ge0,\ k=1,2,\ldots,
\end{equation}
hence \eqnu{currentconservationprf3} is proved for this case.
This completes a proof of \eqnu{currentconservationprf3},
and \eqnu{currentconservation} follows.
\QED
\end{prfof}
Next put
\begin{equation}
\eqna{G}
\begin{array}{l}\dsp
G(\theta)(\gamma,t)=1-\varphi_{\theta}(W,\gamma,t)
\\ \dsp \phantom{G(\theta)(\gamma,t)}
=1-\int_{W\times[y_0,1)}
 \prb{N_{\theta,w,z}(t)=N_{\theta,w,z}(t_0)}\,\mu_0(dw\times dz)
\\ \dsp \phantom{G(\theta)(\gamma,t)}
=y_0+\int_{W\times[0,1)}
 \prb{N_{\theta,w,z}(t)>N_{\theta,w,z}(t_0)}\,\mu_0(dw\times dz)
\\ \dsp
\ \ \gamma=(y_0,t_0)\in\Gamma,\ (\gamma,t)\in\Delta_T\,.
\end{array}
\end{equation}
With \eqnu{Burgersdefdensity} and \eqnu{Poissondecayexplicit},
we have an explicit formula
\begin{equation}
\eqna{HDLop_on_class_of_flows_f}
\begin{array}{l}\dsp
G(\theta)((y_0,0),t)=
1-\int_{W\times[y_0,1)} 
e^{-\Omega_{\theta,w,z}(0,t)}\,\sigma(w,z)\,\lambda(dw)\,dz,
\\ \dsp {}
\gamma=(y_0,0)\in \Gamma_i\,,
\end{array}
\end{equation}
and 
\begin{equation}
\eqna{HDLop_on_class_of_flows_g}
\begin{array}{l}\dsp
G(\theta)((0,t_0),t)=
1-\int_{W\times[0,1)} 
e^{-\Omega_{\theta,w,z}(0,t)}\,\sigma(w,z)\,\lambda(dw)\,dz
\\ \dsp \phantom{G(\theta)((0,t_0),t)=}
-\int_{W\times[0,1)} \sum_{k\ge1}
 \int_{0\le u_1\le \cdots\le u_k\le t_0}
 w(\theta((z,0),u_1),u_1)\, e^{-\Omega_{\theta,w,z}(0,u_1)}\,
\\ \dsp \phantom{G(\theta)((0,t_0),t)=-\int_W}
\times\prod_{i=2}^k \biggl( w(\theta((0,u_{i-1}),u_i),u_i)
\,e^{-\Omega_{\theta,w}(u_{i-1},u_i)}\biggr)\,
\\ \dsp \phantom{G(\theta)((0,t_0),t)=-\int_W}
\times e^{-\Omega_{\theta,w}(u_k,t)}
\prod_{i=1}^kdu_i\, \sigma(w,z)\, \lambda(dw)\,dz,
\\ \dsp \phantom{G(\theta)((0,t_0),t)}
 \gamma=(0,t_0)\in \Gamma_b\cap\Gamma_t\,.
\end{array}
\end{equation}
Using these explicit formula with
\eqnu{Burgersassumps} and \eqnu{Poissonsumrule},
we see from \eqnu{class_of_conti_flows} that $G(\theta)\in\Theta_T$.
Namely, \eqnu{G} defines a map
\begin{equation}
\eqna{Gmap}
G:\ \Theta_T\to\Theta_T
\end{equation}
on the set of flows $\Theta_T$.
\begin{prp}
\prpa{varphi2mu}
If $\theta\in \Theta_T$ is a fixed point of $G$ in \eqnu{Gmap},
namely, if $G(\theta)=\theta$, then 
$\varphi_{\theta}$ defined by \eqnu{varphi_theta} 
uniquely determines
for each $t\in[0,T]$ a probability measure $\mu_{\theta,t}$ on
$W\times [0,1]$ by the equation 
\begin{equation}
\eqna{varphi2mu}
\varphi_{\theta}(dw,\gamma,t)=\mu_{\theta,t}(dw\times[\theta(\gamma,t),1)),
\ \ (\gamma,t)\in \Delta_T\,.
\end{equation}

Furthermore, we have a following formula for 
a change of integration variables:
\begin{equation}
\eqna{varphi2muint}
\begin{array}{l}\dsp
\int_{W\times \{\gamma'\in\Gamma_t\mid \gamma\succeq\gamma'\}}
f(w,\theta(\gamma',t),t) \,\pderiv{\varphi_{\theta}}{\gamma}(dw,\gamma',t)
\,d\gamma'
\\ \dsp
=
\int_{W\times [\theta(\gamma,t),1]} f(w,z,t)\,\mu_{\theta,t}(dw\times dz),
\end{array}
\end{equation}
for integrable function $\dsp f:\ W\times[0,1]\times[0,T]\to\reals$.
\DDD
\end{prp}
\begin{rem}
We are working with a generalization of characteristic curves,
for which $(\gamma,t)$ is a good coordinate.
On the other hand, equation of motions are usually stated in
the space time coordinates $(y,t)$.
This proposition relates the representations in these distinct 
coordinate systems.
\DDD\end{rem}
\begin{prfof}{\protect\prpu{varphi2mu}}
The definition \eqnu{class_of_conti_flows} of $\Theta_T$ implies
that $\dsp\theta(\cdot,t):\ \Gamma_t\to[0,1]$ is continuous, surjective and 
non-decreasing.
Hence to prove \eqnu{varphi2mu}
we only need to prove consistency, namely,
\begin{equation}
\eqna{varphi2muprf1}
\theta(\gamma,t)=\theta(\gamma',t)\ \Rightarrow\ 
\varphi_{\theta}(B,\gamma,t)=\varphi_{\theta}(B,\gamma',t),
\ B\in{\cal B}(W).
\end{equation}
If \eqnu{varphi2muprf1} holds, then \eqnu{varphi2mu} determines
the distribution function 
$\mu_{\theta,t}(B\times[y,1))$, $y\in[0,1]$, on $[0,1]$,
so that $\mu_{\theta,t}(B\times [a,b])$ is determined, and eventually 
$\mu_{\theta,t}$
is determined as a probability measure on the product space 
$W\times[0,1]$.

Assume that $\theta(\gamma,t)=\theta(\gamma',t)$ for
$\gamma=(y_0,t_0)\in\Gamma_t$ and $\gamma'=(y'_0,t'_0)\in\Gamma_t$.
Since $(\Gamma,\succeq)$ is a totally ordered set,
we may assume without loss of generality that $\gamma\succeq\gamma'$.
Then \eqnu{Gammaorder} implies $t_0\le t'_0\le t$ and $y_0\ge y'_0$.
Non-decreasing property of the point process $N_{\theta,w,z}(t)$ in $t$ and
monotonicity of measures imply with the definition \eqnu{varphi_theta}, 
\begin{equation}
\eqna{varphi2muprf2}
\varphi_{\theta}(B,\gamma',t)\ge \varphi_{\theta}(B,\gamma,t),
\ \mbox{ and }\ 
\varphi_{\theta}(B^c,\gamma',t)\ge \varphi_{\theta}(B^c,\gamma,t),
\end{equation}
for any $B\in {\cal B}(W)$,

On the other hand, the assumption $G(\theta)=\theta$ and \eqnu{G}
and $\theta(\gamma,t)=\theta(\gamma',t)$ imply
\[ \begin{array}{l}\dsp
\varphi_{\theta}(B,\gamma,t)+\varphi_{\theta}(B^c,\gamma,t)
=\varphi_{\theta}(W,\gamma,t)=1-G(\theta)(\gamma,t)=1-\theta(\gamma,t)
\\ \dsp
= 1-\theta(\gamma',t)=
\varphi_{\theta}(B,\gamma',t)+\varphi_{\theta}(B^c,\gamma',t).
\end{array} \]
Hence
\begin{equation}
\eqna{varphi2muprf3}
 \varphi_{\theta}(B,\gamma',t)-\varphi_{\theta}(B,\gamma,t)
=-( \varphi_{\theta}(B^c,\gamma',t)-\varphi_{\theta}(B^c,\gamma,t) ).
\end{equation}
Combining
\eqnu{varphi2muprf2} and
\eqnu{varphi2muprf3}, we see that \eqnu{varphi2muprf1} holds,
which implies \eqnu{varphi2mu}.

Next to prove \eqnu{varphi2muint}, we first prove
\begin{equation}
\eqna{fundamentallawofanalysisforvarphi}
\int_{\gamma\succeq\gamma'} \pderiv{\varphi_{\theta}}{\gamma}(B,\gamma',t)
\,d\gamma' 
= \mu_{\theta,t}(B\times[\theta(\gamma,t),1)),
\ \ (\gamma,t)\in\Delta_T\,,\ B\in{\cal B}(W).
\end{equation}
In fact, note that \eqnu{varphi_theta} implies
\begin{equation}
\eqna{varphi_zero}
\varphi_{\theta}(B,(1,0),t)=0,\ \ B\in{\cal B}(W),\ t\ge0.
\end{equation}
Then, if $\gamma=(y_0,0)\in\Gamma_i$\,,
The explicit formula \eqnu{Poissondecayexplicit},
with the definitions
\eqnu{varphi_theta}
\eqnu{dphidgamma}
\eqnu{Gammaintegral},
implies
\[ \begin{array}{l}\dsp
\int_{\gamma\succeq\gamma'} \pderiv{\varphi_{\theta}}{\gamma}(B,\gamma',t)
\,d\gamma' = 
-\int_{y_0}^1 \pderiv{\varphi_{\theta}}{z}(B,(z,0),t)\,dz
\\ \dsp \phantom{%
\int_{\gamma\succeq\gamma'} \pderiv{\varphi_{\theta}}{\gamma}(B,\gamma',t)
\,d\gamma'}
=\varphi(B,(y_0,0),t),
\end{array} \]
where we used \eqnu{varphi_zero} in the last line.
If $\gamma=(0,t_0)\in\Gamma_b$\,,
the explicit formula \eqnu{Poissondecayexplicit}
similarly implies
\[ \begin{array}{l}\dsp
\int_{\gamma\succeq\gamma'} \pderiv{\varphi_{\theta}}{\gamma}(B,\gamma',t)
\,d\gamma' 
\\ \dsp \phantom{}
= -\int_0^1 \pderiv{\varphi_{\theta}}{z}(B,(z,0),t)\,dz
+\int_0^{t_0}  \pderiv{\varphi_{\theta}}{u}(B,(0,u),t)\,du
\\ \dsp \phantom{}
=-(\varphi_{\theta}(B,(1,0),t)-\varphi_{\theta}(B,(0,0),t))
+(\varphi_{\theta}(B,(0,t_0),t)-\varphi_{\theta}(B,(0,0),t))
\\ \dsp \phantom{}
=\varphi(B,(0,t_0),t),
\end{array}\]
hence \eqnu{fundamentallawofanalysisforvarphi} follows from
\eqnu{varphi2mu}.

Define a measure $\nu_{\theta,t}$ on $W\times \Gamma_t$ by
\[
\nu_{\theta,t}(dw\times d\gamma)
=\pderiv{\varphi_{\theta}}{\gamma}(dw,\gamma,t)\,d\gamma.
\]
Then \eqnu{fundamentallawofanalysisforvarphi} implies
\[
\nu_{\theta,t}(B\times \{\gamma'\mid \gamma\succeq\gamma'\})
= \mu_{\theta,t}(B\times[\theta(\gamma,t),1)),
\ \gamma\in\Gamma_t\,.
\]
This implies that,
if we put
$\dsp X_{\theta,t}=({\rm id}_W,\theta)$,
where ${\rm id}_W$ is the identity map on $W$,
then 
$\mu_{\theta,t}$ is the image measure of $\nu_{\theta,t}$ 
with respect to the map 
$X_{\theta,t}:\ W\times\Gamma_t\to W\times[0,1]$: 
\[\mu_{\theta,t}=\nu_{\theta,t}\circ X^{-1}_{\theta,t}\,. \]
Therefore \eqnu{varphi2muint} follows.
\QED
\end{prfof}
\begin{thm}
\thma{soln}
Assume that $G$ in \eqnu{Gmap} has a fixed point,
and denote the fixed point by $\dsp y_C\in\Theta_T$;
\begin{equation}
\eqna{yCFP}
y_C=G(y_C)
\end{equation}
Put $\dsp\varphi=\varphi_{y_C}$ and $\dsp \mu_t=\mu_{y_C,t}$,
where $\dsp\varphi_{\theta}$ and $\dsp\mu_{\theta,t}$ are defined
by \eqnu{varphi_theta} and \eqnu{varphi2mu} with $\dsp\theta=y_C$\,.
Then the so defined $y_C$ and $\mu_t$ satisfy all the properties stated in
\thmu{Burgers}.
\DDD
\end{thm}
\begin{rem}
With \eqnu{varphi_theta} and \eqnu{varphi2mu}, the theorem implies 
an expression 
\begin{equation}
\eqna{explicit_soln}
\begin{array}{l}\dsp
\mu_{y_C,t}(dw\times[y_C(\gamma,t),1))
=\int_{z\in[y_0,1)}
 \prb{N_{y_C,w,z}(t)=N_{y_C,w,z}(t_0)}\,\mu_0(dw\times dz),
\\ \dsp
\gamma=(y_0,t_0)\in\Gamma_t\,,
\end{array}
\end{equation}
for the solution to \thmu{Burgers}.

The properties of the solution claimed in \thmu{Burgers} are mostly
contained in the previous propositions and explicit formulas.
The remaining point is the Lipschitz continuity of $\dsp\mu_t$:
If $w(y,t)=0$ for certain time interval for all $w\in W$,
then the characteristic curve $y_C$ will remain constant for the interval,
and a small change in $y=y_C(\gamma,t)$ may result in a large change
in $\gamma$ and hence in $\mu_t$\,.
It turns out that the situation causes no problem, because then
the change in the quantity $\varphi$ and eventually $\mu_t$ are small
in the time interval, hence continuity follows.
\DDD\end{rem}
\begin{prfof}{\protect\thmu{soln}}
By definition \eqnu{G}, we have
\begin{equation}
\eqna{yCvarphi}
y_C(\gamma,t)=G(y_C)(\gamma,t)=1-\varphi_{y_C}(W,\gamma,t),
\end{equation}
for which monotonicity properties stated in \thmu{Burgers}
are direct consequences of \eqnu{varphi_theta}.
\prpu{finiteevaporation} and \prpu{currentconservation} with
explicit formula \eqnu{HDLop_on_class_of_flows_f} and 
\eqnu{HDLop_on_class_of_flows_g}, or \eqnu{varphi_theta} and
\eqnu{dphidu},
imply that $y_C$ and $\dsp \pderiv{y_C}{t}$ is continuous in $(\gamma,t)$.

Let $h:\ W\to\reals$ be a measurable function
satisfying \eqnu{Burgersprf3}.
For $\dsp (y,t),\; (y',t)\in [0,1]\times[0,T]$,
choose $\gamma\in\Gamma_t$ and $\gamma'\in\Gamma_{t}$ such that
$\dsp y=y_C(\gamma,t)$ and $\dsp y'=y_C(\gamma',t)$.
We may assume $\gamma\succeq\gamma'$.
Then monotonicity of $y_C$ implies $y\le y'$, 
hence $h(w)\le 1$ and \eqnu{soliditycond} imply
\[ \begin{array}{l}\dsp
\int_W h(w)\mu_{t}(dw\times[y,1))-\int_W h(w)\mu_{t}(dw\times[y',1))
=
\int_W h(w)\mu_{t}(dw\times[y,y'))
\\ \dsp {}
\le \mu_t(W\times [y,y'))=y'-y.
\end{array} \]
Next,
for $\dsp (y,t),\; (y,t')\in [0,1]\times[0,T]$,
choose $\gamma\in\Gamma_t$ and $\gamma'\in\Gamma_{t'}$ such that
$\dsp y=y_C(\gamma,t)$ and $\dsp y=y_C(\gamma',t')$.
Then, since,
by definition \eqnu{varphi_theta}, $\varphi(B,\gamma,t)$ is monotone
also in $t$,
\[ \begin{array}{l}\dsp
\left|\int_W h(w)\mu_{t}(dw\times[y,1))-\int_W h(w)\mu_{t'}(dw\times[y,1))
\right|
\\ \dsp {}
=\biggl|\int_W h(w)\,(\varphi(dw,\gamma,t)-\varphi(dw,\gamma',t'))\biggr|
\\ \dsp {}
\le |\varphi(W,\gamma,t)-\varphi(W,\gamma',t)|
+ |\varphi(W,\gamma',t)-\varphi(W,\gamma',t')|
\\ \dsp {}
=|y_C(\gamma,t)-y_C(\gamma',t)|
+|y_C(\gamma',t)-y_C(\gamma',t')|
\\ \dsp {}
=2 |y-y_C(\gamma',t)|
\\ \dsp {}
=2 |y_C(\gamma',t')-y_C(\gamma',t)|
\\ \dsp {}
\le 2\sup_{(\gamma'',t'')\in\Delta_T} 
\biggl|\pderiv{y_C}{t}(\gamma'',t'')\biggr|\,|t'-t|.
\end{array}\]
Using \prpu{currentconservation}, monotonicity, 
and \prpu{finiteevaporation}, 
with \eqnu{yCvarphi}, we have
\[ \begin{array}{l}\dsp
\pderiv{y_C}{t}(\gamma,t)
= \int_{\gamma\succeq\gamma'} \int_W w(y_C(\gamma',t),t)\,
\pderiv{\varphi}{\gamma}(dw,\gamma',t')\,d\gamma'
\\ \dsp{}
\le \int_{\gamma'\in\Gamma_t} \int_W \wnorm{w}
\pderiv{\varphi}{\gamma}(dw,\gamma',t)\,d\gamma'
\\ \dsp{}
\le M_W e^{2C_Wt},
\end{array}\]
which proves \eqnu{globalLipschitz}.

The initial conditions \eqnu{U0} and surjectivity 
for $y_C\in\Theta_T$ is contained
in the definition \eqnu{class_of_conti_flows} of $\Theta_T$,
and the initial condition \eqnu{Burgersdefdensity}, for 
for $\mu_t$ is in the definition \eqnu{varphi_theta} with \prpu{varphi2mu}.

For $t>0$ and $B\in{\cal B}(W)$,
\eqnu{varphi2mu}, \eqnu{varphi_theta}, \eqnu{Burgersdefdensity}, and
\eqnu{Burgersdensity1} imply
\[
\mu_t(dw\times[0,1))=\varphi(dw,(0,t),t)
=\mu_0(dw\times [0,1))
=\lambda(dw),
\]
which proves \eqnu{conservation}.

For $(y,t)\in[0,1]\times[0,T]$ choose $\gamma\in\Gamma_t$
such that $\dsp y=y_C(\gamma,t)$. 
Then \eqnu{varphi2mu} and \eqnu{G} imply
\[
\mu_t(W\times[y,1))=\varphi(W,\gamma,t)
=1-G(y_C)(\gamma,t)=1-y_C(\gamma,t)=1-y
\]
which proves \eqnu{soliditycond}.

To prove evolution equation \eqnu{phic},
\prpu{varphi2mu} implies, for $\gamma=(y_0,t_0)\in\Gamma_t$\,,
\[
\mu_{t}(dw\times[y_C(\gamma,t),1))-\mu_{t_0}(dw\times[y_0,1))
=\varphi(dw,\gamma,t)-\varphi(dw,\gamma,t_0),
\]
for which \prpu{currentconservation} and \eqnu{varphi2muint} further imply
\[ \begin{array}{l}\dsp
\mu_{t}(dw\times[y_C(\gamma,t),1))-\mu_{t_0}(dw\times[y_0,1))
\\ \dsp {}
=-\int_{t_0}^t\biggl(\int_{\gamma\succeq\gamma'}
w(y_C(\gamma',s),s)\,\pderiv{\varphi}{\gamma}(dw,\gamma',s)\,d\gamma'\biggr)
\,ds
\\ \dsp {}
=-\int_{t_0}^t\biggl(\int_{y_C(\gamma,s)}^1w(z,s)\, \mu_s(dw\times dz)\biggr)
\,ds,
\end{array}\]
which proves \eqnu{phic}.
\QED
\end{prfof}

\section{Fixed point and existence of solution.}
\seca{FP}

In \thmu{soln} we assumed existence of a fixed point $\theta=y_C$
of a map $G$ defined in \eqnu{G}.
To complete a proof of existence of a solution for \thmu{Burgers},
we prove that \eqnu{G} has a fixed point.
In fact, the assumptions \eqnu{rw} and \eqnu{Cw} on $W$ imply that
that the fixed point is unique.
This is the core of the existence proof for \thmu{Burgers}, 
and we heavily rely on the explicit formulas
\eqnu{HDLop_on_class_of_flows_f} and \eqnu{HDLop_on_class_of_flows_g}.
\begin{thm}
\thma{FPtheorem}
The map $G:\ \Theta\to\Theta$ in \eqnu{Gmap} has a unique fixed point
$\dsp y_C\in\Theta_T$\,, namely, there is a unique $y_C$ 
which satisfies \eqnu{yCFP}.
\DDD
\end{thm}
\begin{prf}
For $t\ge0$ and $\theta$ and $\theta'$ in $\Theta_T$ define
\begin{equation}
\eqna{theta_norm}
d(\theta',\theta,t)=\sup_{\gamma\in\Gamma_t} 
|\theta'(\gamma,t)-\theta(\gamma,t)|.
\end{equation}
We first accumulate basic formulas for evaluating
$\omega_{\theta,w,z}$ in \eqnu{w2omega}
and $\Omega_{\theta,w,z}$ in \eqnu{omegaOmega}.
In the following lemma, we write $\Omega_{\theta,w,z}(s,t)$
also for $s>0$ in \eqnu{omegaOmega} whenever it becomes notationally
simpler, though the quantity is actually independent of $z$ for $s>0$.
Recall the notation $\tilde{\Omega}_w(s,t)$ in \eqnu{wstar}.
\begin{lem}
\lema{FPtheoremprf1}
Let $\theta\in\Theta_T$ and $\theta'\in\Theta_T$\,.
Then for $\dsp (\gamma,t)\in\Delta_T$ with $\gamma=(z,s)$, we have
\begin{enumerate}
\item
$\dsp w(1,t)-C_W\le w(\theta(\gamma,t),t) \le w(1,t)+C_W$,

\item
$\dsp 0<e^{-\Omega_{\theta,w,z}(s,t)}
\le e^{-\tilde{\Omega}_w(s,t)+C_W\,(t-s)}$,

\item
$\dsp |w(\theta'(\gamma,t),t)-w(\theta(\gamma,t),t)|
\le C_W\,d(\theta',\theta,t)$,

\item
$\dsp |e^{-\Omega_{\theta',w,z}(s,t)}-e^{-\Omega_{\theta,w,z}(s,t)}|
\le C_W\,e^{-\tilde{\Omega}_w(s,t)+C_W\,(t-s)}
\int_s^t d(\theta',\theta,v)\,dv$.

\end{enumerate}
\DDD
\end{lem}
\begin{prf}
The first estimate is an elementary consequence of \eqnu{Cw}
and the mean value theorem, if one notes that $\theta(\gamma,t)\in[0,1]$.
This and the definitions \eqnu{omegaOmega} and \eqnu{wstar},
and non-negativity of $w\in W$ leads to the second estimate.
With the definition \eqnu{theta_norm},
the third estimate is similarly proved as the first one.
The last estimate follows from these estimates and
\[
|e^{-x'}-e^{-x}|=|e^{-(x'\vee x)}-e^{-(x'\wedge x)}|
= e^{-(x'\wedge x)} \, (1- e^{-|x'-x|})
\\ \dsp \le e^{-(x'\wedge x)} |x'-x|.
\]
\QED
\end{prf}
\begin{lem}
\lema{FPtheoremprf2}
It holds that
\begin{equation}
\eqna{FPtheoremprf2}
d(G(\theta'),G(\theta),t)\le 2C_W\,e^{2C_W\,T}
\int_0^t d(\theta',\theta,v)\,dv,
\ \ \theta,\,\theta\in\Theta_T\,,\ t\in[0,T].
\end{equation}
\DDD
\end{lem}
\begin{prf}
If $\gamma=(y_0,0)\in\Gamma_i$, then
applying \lemu{FPtheoremprf1} to \eqnu{HDLop_on_class_of_flows_f},
we have
\[\begin{array}{l}\dsp
|G(\theta')(\gamma,t)-G(\theta)(\gamma,t)|
\\ \dsp \phantom{}
\le C_W\,e^{C_W\,t} \int_{W\times[y_0,1]}e^{-\tilde{\Omega}_w(0,t)}
\sigma(w,z)\lambda(dw)\,dz \, \int_0^td(\theta',\theta,v)\,dv,
\ \ \gamma\in\Gamma_i\,\ t\in[0,T].
\end{array}\]
Non-negativity of $\tilde{\Omega}_w$ and \eqnu{Burgersassumps} and
the fact that $\lambda$ is a probability measure further leads to
\begin{equation}
\eqna{FPtheoremprf3}
\sup_{\gamma\in\Gamma_i}|G(\theta')(\gamma,t)-G(\theta)(\gamma,t)|
\le C_W\,e^{C_W\,T} \, \int_0^td(\theta',\theta,v)\,dv,
\ \ t\in[0,T].
\end{equation}
The rest of the proof is for the case 
$\gamma=(0,t_0)\in\Gamma_t\cap\Gamma_b$\,.
On applying \lemu{FPtheoremprf1} to each term of 
\eqnu{HDLop_on_class_of_flows_g}, we use an elementary equality
\begin{equation}
\eqna{FPtheoremprf4}
\prod_{i=1}^n b_i -\prod_{i=1}^na_i= \sum_{j=1}^n 
\biggl(\prod_{i=1}^{j-1}b_i\biggr)\,(b_j-a_j)\,
\biggl(\prod_{i=j+1}^{n}a_i\biggr),
\end{equation}
where, and in the following, we adopt a notation
\[
 \prod_{i=1}^0 b_i=\prod_{i=n+1}^n a_i=1
\]
to simplify the formulas.
We apply \eqnu{FPtheoremprf4} to the difference
of \eqnu{HDLop_on_class_of_flows_g} and its analog, with 
$\theta$ replaced by $\theta'$, where
$b_i$'s are the factors depending on $\theta'$,
and $a_i$'s the factors depending on $\theta$.
We then apply the last $2$ estimates in \lemu{FPtheoremprf1}
to the factor of the form $b_i-a_i$, and 
apply the first $2$ estimates to other factors.
We have
\begin{equation}
\eqna{FPtheoremprf5}
|G(\theta')((0,t_0),t)-G(\theta)((0,t_0),t)|\le I_1(t_0,t)+I_2(t_0,t),
\end{equation}
where
\begin{equation}
\eqna{FPtheoremprf6}
\begin{array}{l}\dsp
I_1(t_0,t)
=C_We^{C_W\,t}\int_W e^{-\tilde{\Omega}_w(0,t)}
\,\biggl(1
\\ \dsp \phantom{I_1(t_0,t)=C_We^{C_W\,t}\int_W}
+\sum_{k\ge1} \int_{0\le u_1\le\ldots\le u_k\le t_0}
\prod_{i=1}^k (w(1,u_i)+C_W)\, \prod_{i=1}^kdu_i\biggr)
\lambda(dw)
\\ \dsp \phantom{I_1(t_0,t)=}
\times\int_0^td(\theta',\theta,v)\,dv,
\end{array}
\end{equation}
and
\begin{equation}
\eqna{FPtheoremprf7}
\begin{array}{l}\dsp
I_2(t_0,t)=C_We^{C_W\,t}\int_We^{-\tilde{\Omega}_w(0,t)}
\times\sum_{k\ge1} \int_{0\le u_1\le\ldots\le u_k\le t_0}
\sum_{j=1}^k\biggl(
 d(\theta',\theta,u_j)
\\ \dsp \phantom{=I_2(t_0,t)=C_We^{C_W\,t}\int_W}
\times\prod_{i;\;1\le i \le k,\,i\ne j}
(w(1,u_i)+C_W)\biggr)\,\prod_{i=1}^k du_i
\lambda(dw).
\end{array}
\end{equation}

We apply \eqnu{integralopcompactness}
to \eqnu{FPtheoremprf6}, to find
\begin{equation}
\eqna{FPtheoremprf8}
\begin{array}{l}\dsp
I_1(t_0,t)=
C_We^{C_W\,t}\int_We^{-\tilde{\Omega}_w(0,t)}
e^{\tilde{\Omega}_w(0,t_0)+C_W\,t_0}\lambda(dw)
\times\int_0^td(\theta',\theta,v)\,dv
\\ \dsp \phantom{I_1(t_0,t)}
=
C_We^{C_W\,(t+t_0)}\int_We^{-\tilde{\Omega}_w(t_0,t)}\lambda(dw)
\,\int_0^td(\theta',\theta,v)\,dv,
\\ \dsp \phantom{I_1(t_0,t)}
\le
C_We^{2C_W\,T}\,\int_0^td(\theta',\theta,v)\,dv.
\end{array}
\end{equation}
To evaluate \eqnu{FPtheoremprf7}, we first change
an integration variable $u_j$ to $v$ and
change the order of summation for $j$ and $k$, to find
\[ \begin{array}{l}\dsp
I_2(t_0,t)
\\ \dsp \phantom{}
=C_We^{C_W\,t}\int_We^{-\tilde{\Omega}_w(0,t)}
\int_0^{t_0}d(\theta',\theta,v)
\\ \dsp \phantom{=}
\times\biggl(\sum_{j\ge1}
 \int_{0\le u_1\le\ldots\le u_{j-1}\le v}
\prod_{i=1}^{j-1}(w(1,u_i)+C_W)\prod_{i=1}^{j-1} du_i\biggr)
\\ \dsp \phantom{=}
\times\biggl(\sum_{k\ge j}
 \int_{v\le u_{j+1}\le \ldots \le u_k\le t_0}
\prod_{i=j+1}^{k}(w(1,u_i)+C_W)\prod_{i=j+1}^k du_i\biggr)
\,\lambda(dw)\, dv.
\end{array} \]
We apply \eqnu{integralopcompactness}
to the summation in $j$ and to the summation in $k$, to find
\begin{equation}
\eqna{FPtheoremprf9}
\begin{array}{l}\dsp
I_2(t_0,t)
\\ \dsp \phantom{}
=C_We^{C_W\,t}\int_We^{-\tilde{\Omega}_w(0,t)}
\int_0^{t_0}d(\theta',\theta,v)
e^{\tilde{\Omega}_w(0,v)+C_W\,v}\,e^{\tilde{\Omega}_w(v,t_0)+C_W\,(t_0-v)}
\,\lambda(dw)\, dv
\\ \dsp \phantom{}
=C_We^{C_W\,(t+t_0)}\int_We^{-\tilde{\Omega}_w(t_0,t)}\,\lambda(dw)
\times \int_0^{t_0}d(\theta',\theta,v)dv
\\ \dsp \phantom{}
\le C_We^{2C_W\,T} \int_0^{t}d(\theta',\theta,v)dv.
\end{array}
\end{equation}
The equations 
\eqnu{FPtheoremprf5}, \eqnu{FPtheoremprf8}, and \eqnu{FPtheoremprf9}
imply
\begin{equation}
\eqna{FPtheoremprf10}
\sup_{\gamma\in\Gamma_b\cap\Gamma_t}
|G(\theta')((0,t_0),t)-G(\theta)((0,t_0),t)|\le 
2C_We^{2C_W\,T} \int_0^{t}d(\theta',\theta,v)dv.
\end{equation}
The equations 
\eqnu{theta_norm},
\eqnu{FPtheoremprf3}, and \eqnu{FPtheoremprf10}
finally imply \eqnu{FPtheoremprf2}.
\QED
\end{prf}
Let us continue the proof of \thmu{FPtheorem}.

Define $\theta_0\in\Theta_T$ by
\begin{equation}
\eqna{theta0}
\theta_0((y_0,t_0),t)=y_0\,,\ \ ((y_0,t_0),t)\in\Delta_T\,,
\end{equation}
and define a sequence of flows $\theta_k\in\Theta_T$, $k\in\pintegers$,
inductively by \eqnu{theta0} and
\begin{equation}
\eqna{thetak}
\theta_{k+1}=G(\theta_k),\ \ k\in\pintegers.
\end{equation}
\lemu{FPtheoremprf2} implies
\begin{equation}
\eqna{FPtheoremprf11}
\begin{array}{l}\dsp
d(\theta_{k+1},\theta_k,t)=d(G(\theta_k),G(\theta_{k-1}),t)
\le C\,\int_0^t d(\theta_k,\theta_{k-1},v)dv,
\\ \dsp \phantom{}
t\in[0,T],\ k=1,2,\ldots.
\end{array}
\end{equation}
where $\dsp C=2C_W\,e^{2C_W\,T}$.
Iterating, we obtain estimates which, by induction, 
is seen to have an expression
\begin{equation}
\eqna{FPtheoremprf12}
\begin{array}{l}\dsp
d(\theta_{k+1},\theta_k,t)
\\ \dsp \phantom{}
\le C^k\int_{0\le u_1\le \ldots \le u_k\le t}d(\theta_1,\theta_0,u_1)\,
\prod_{i=1}^k du_i
= C^k\int_0^t \frac{(t-u)^{k-1}}{(k-1)!}\,d(\theta_1,\theta_0,u)\,du,
\\ \dsp \phantom{}
t\in[0,T],\ k\in\pintegers.
\end{array}
\end{equation}
Since $\theta\in\Theta_T$ takes values in $[0,1]$,
the definition \eqnu{theta_norm} implies
\begin{equation}
\eqna{FPtheoremprf13}
d(\theta',\theta,t)\le1,\ \ \theta,\theta'\in\Theta_T\,,\ t\in[0,T].
\end{equation}
Substituting \eqnu{FPtheoremprf13} in \eqnu{FPtheoremprf12},
\begin{equation}
\eqna{FPtheoremprf14}
d(\theta_{k+1},\theta_k,t)\le \frac{(C\,t)^k}{k!}\,,\ k\in\pintegers.
\end{equation}
Since the summation in $k$ of the right hand side of 
\eqnu{FPtheoremprf14} converges to $\dsp e^{C\,t}$,
\[\theta_k(\gamma,t)=\theta_0(\gamma,t)+\sum_{i=0}^{k-1}
 (\theta_{i+1}(\gamma,t)-\theta_i(\gamma,t)) \]
converges uniformly in $(\gamma,t)\in\Delta_T$\,.
Denote the limit as
\begin{equation}
\eqna{FPtheoremprf15}
y_C(\gamma,t)=\limf{k} \theta_k(\gamma,t),\ \ (\gamma,t)\in\Delta_T\,.
\end{equation}
The equations \eqnu{FPtheoremprf15}, \eqnu{thetak}, and \eqnu{FPtheoremprf14}
imply \eqnu{yCFP}.
Since $\theta_k\in\Theta_T$ for all $k$,
$y_C$ also takes values in $[0,1]$,
non-decreasing in $\gamma$ for each $t$, and
non-decreasing in $t$ for each $\gamma$.
Since the convergence \eqnu{FPtheoremprf15} is uniform in $(\gamma,t)$,
and $\theta_k\in\Theta_T$ are continuous, $y_C$ is also continuous.
Also 
\[ y_C((y_0,t_0),t_0)=y_0\,,\ \ (y_0,t_0)\in\Gamma_T\,, \]
holds.
In particular,
$\dsp y_C((0,t),t)=0$ holds, and also \eqnu{Gamma_infinity} implies
$\dsp y_C((1,0),t)=1$, hence with continuity, $y_C$ is surjective
in $\gamma$ for each $t$.
This proves $y_C\in\Theta_T$\,, namely, existence of a fixed point
of $G$ in $\Theta_T$\,.

Suppose there is another fixed point $\tilde{y}_C\in\Theta_T$ of $G$.
Then \eqnu{yCFP} and \lemu{FPtheoremprf2} imply
\[
d(\theta',\theta,t)=d(G(\theta'),G(\theta),t)
\le C\,\int_0^t d(\theta',\theta,v)\,dv,\ \ t\in[0,T].
\]
where $\dsp C=2C_W\,e^{2C_W\,T}$.
Gronwall's inequality implies $d(\theta',\theta,t)=0$, $t\in[0,T]$.
Namely, $\theta'=\theta$.
This proves uniqueness of the fixed point of $G$.
\QED
\end{prf}

\section{Uniqueness of the solution.}
\seca{prfofBurgers}

In previous sections we proved existence of a solution
$(y_C,\mu_t)$ in \thmu{Burgers}.
In this section we complete the proof of \thmu{Burgers} by
proving that the solution is unique.
Assume that $(y_C,\mu_t)$ and $(\tilde{y}_C,\tilde{\mu}_t)$ are
the pairs which satisfy all the properties stated in \thmu{Burgers}.

Fix $\gamma=(y_0,t_0)\in\Gamma$\,, and let $t\ge t_0$\,.
Since by assumption $(\mu_t,y_C)$
satisfy equation of motion \eqnu{phic} with initial and boundary conditions
\eqnu{U0} and \eqnu{conservation},
\begin{equation}
\eqna{Burgersprf1}
\begin{array}{l}\dsp
\mu_t(dw\times[y_C(\gamma,t),1))
\\ \dsp
=\mu_{t_0}(dw\times[y_0,1))
\\ \dsp \phantom{=}
+\int_{t_0}^t \int_{z\in[y_C(\gamma,s),1)} (w(1,s)-w(z,s))
 \mu_s(dw\times dz)\,ds
\\ \dsp \phantom{=}
-\int_{t_0}^t w(1,s)
 \mu_s(dw\times [y_C(\gamma,s),1))\,ds,\ \ t\ge t_0\,.
\end{array}
\end{equation}
Note that \eqnu{conservation} and \eqnu{rw} do not rule out
a possibility that $\mu_t$ has an unbounded support concerning $\wnorm{w}$.
Therefore, a direct application of Gronwall type inequalities to 
the last term in the right hand side of \eqnu{Burgersprf1}
may lead to divergent expression upon integration with respect to $w$.
We work around this problem by the following.
\begin{lem}
\lema{Burgersprf}
\begin{equation}
\eqna{Burgersprf10}
\begin{array}{l}\dsp
\mu_t(dw\times[y_C(\gamma,t),1))
=e^{-\tilde{\Omega}_w(t_0,t)}\mu_{t_0}(dw\times[y_0,1))
\\ \dsp \phantom{\biggl|}
+\int_{t_0}^t 
e^{-\tilde{\Omega}_w(s,t)}\int_{x\in[y_C(\gamma,s),1)}  \pderiv{w}{z}(x,s)
 \mu_s(dw\times [y_C(\gamma,s),x))\,dx\,ds,
\end{array}
\end{equation}
where $\tilde{\Omega}_w(s,t)$ is as in \eqnu{wstar}.
\DDD
\end{lem}
\begin{prf}
Iterating \eqnu{Burgersprf1} and using Fubini's Theorem, we have
\begin{equation}
\eqna{Burgersprf2}
\begin{array}{l}\dsp
\mu_t(dw\times[y_C(\gamma,t),1))
\\ \dsp
=\mu_{t_0}(dw\times[y_0,1))\,
\sum_{\ell=0}^k\frac1{\ell!}\,(-\tilde{\Omega}_w(t_0,t))^{\ell}
\\ \dsp \phantom{=}
+\int_{t_0}^t 
\sum_{\ell=0}^k\frac1{\ell!}\,(-\tilde{\Omega}_w(s,t))^{\ell}
\int_{z\in[y_C(\gamma,s),1)} (w(1,s)-w(z,s))
 \mu_s(dw\times dz)\,ds
\\ \dsp \phantom{=}
-\int_{t_0}^t w(1,s)
\frac1{k!}\,(-\tilde{\Omega}_w(s,t))^k
 \mu_s(dw\times [y_C(\gamma,s),1))\,ds,
\\ \dsp \phantom{}
t\ge t_0\,,\ k=0,1,2,\ldots.
\end{array}
\end{equation}
Since $w\in W$ are non-negative valued, so are $\tilde{\Omega}_w(s,t)$ for $s\le t$
and
\begin{equation}
\eqna{Burgersprf4}
0\le \tilde{\Omega}_w(s,t)\le \wnorm{w}\,T,\ \ 0\le s\le t\le T.
\end{equation}
Taylor's Theorem and \eqnu{Burgersprf4} imply
\begin{equation}
\eqna{Burgersprf5}
\begin{array}{l}\dsp
\biggl|e^{-\tilde{\Omega}_w(s,t)}
-\sum_{\ell=0}^k \frac1{\ell!}(\tilde{\Omega}_w(s,t))^{\ell}\biggr|
\le \frac1{(k+1)!}(\wnorm{w}\,T)^{k+1},
\\ \dsp 
0\le s\le t\le T,\ \ k=0,1,2,\ldots.
\end{array}
\end{equation}
Note also that
\eqnu{Cw} and the mean value theorem imply
\begin{equation}
\eqna{Burgersprf6}
|w(1,s)-w(z,s)|\le C_W\,, \ \ z\in[0,1],\ s\in[0,T].
\end{equation}

Fix a constant $M>T$ (e.g., $M=2T$), and let $B\in{\cal B}(W)$.
Then, monotonicity of measures,
\eqnu{Burgersprf2},
\eqnu{Burgersprf4},
\eqnu{Burgersprf5},
\eqnu{Burgersprf6},
and \eqnu{conservation}
imply
\begin{equation}
\eqna{Burgersprf7}
\begin{array}{l}\dsp
\biggl|\int_B e^{-M\wnorm{w}}\mu_t(dw\times[y_C(\gamma,t),1))
-\int_B e^{-M\wnorm{w}}\,e^{-\tilde{\Omega}_w(t_0,t)}
\mu_{t_0}(dw\times[y_0,1))
\\ \dsp \phantom{\biggl|}
-\int_{t_0}^t 
\int_B e^{-M\wnorm{w}}\,e^{-\tilde{\Omega}_w(s,t)}
\int_{z\in[y_C(\gamma,s),1)} (w(1,s)-w(z,s))
\mu_s(dw\times dz)\,ds
\biggr|
\\ \dsp\phantom{}
\le\frac{T^{k+1}}{(k+1)!}\,(1+C_W\,T+k+1)
\int_B e^{-M\wnorm{w}} \wnorm{w}^{k+1}\,\lambda(dw)
\\ \dsp\phantom{}
\le\frac{T^{k+1}}{(k+1)!}\,(C_W\,T+k+2)
\,\lambda(B)\,\sup_{x\ge0} x^{k+1}\,e^{-M\,x}
\\ \dsp\phantom{}
\le\frac{T^{k+1}}{(k+1)!}\,(C_W\,T+k+2)
\,\sup_{x\ge0} x^{k+1}\,e^{-M\,x},
\\ \dsp\phantom{}
t\ge t_0\,,\ k=0,1,2,\ldots.
\end{array}
\end{equation}
By elementary calculus,
\begin{equation}
\eqna{Burgersprf8}
\begin{array}{l}\dsp
\log \sup_{x\ge0} x^{k+1}\,e^{-M\,x}
=(k+1)\,(\log\frac{k+1}M-1)
\\ \dsp \phantom{}
= \int_0^{k+1}\log y\,dy -(k+1)\log M
\\ \dsp \phantom{}
\le \sum_{\ell=1}^{k+1}\log \ell -(k+1)\log M
=\log\frac{(k+1)!}{M^{k+1}}\,.
\end{array}
\end{equation}
Combining \eqnu{Burgersprf7} and \eqnu{Burgersprf8}, we have
\[\begin{array}{l}\dsp
\biggl|\int_B e^{-M\wnorm{w}}\mu_t(dw\times[y_C(\gamma,t),1))
-\int_B e^{-M\wnorm{w}}\,e^{-\tilde{\Omega}_w(t_0,t)}
\mu_{t_0}(dw\times[y_0,1))
\\ \dsp \phantom{\biggl|}
-\int_{t_0}^t 
\int_B e^{-M\wnorm{w}}\,e^{-\tilde{\Omega}_w(s,t)}
\int_{z\in[y_C(\gamma,s),1)} (w(1,s)-w(z,s))
 \mu_s(dw\times dz)\,ds
\biggr|
\\ \dsp\phantom{}
\le \biggl(\frac{T}{M}\biggr)^{k+1}(C_WT+k+2),
\ \ t\ge t_0\,,\ B\in{\cal B}(W),\ k=0,1,2,\ldots,
\end{array}\]
which implies, by fixing $M>T$ and considering $k\to\infty$,
\[\begin{array}{l}\dsp
\int_B e^{-M\wnorm{w}}\mu_t(dw\times[y_C(\gamma,t),1))
=\int_B e^{-M\wnorm{w}}\,e^{-\tilde{\Omega}_w(t_0,t)}
\mu_{t_0}(dw\times[y_0,1))
\\ \dsp \phantom{\biggl|}
+\int_{t_0}^t 
\int_B e^{-M\wnorm{w}}\,e^{-\tilde{\Omega}_w(s,t)}
\int_{z\in[y_C(\gamma,s),1)} (w(1,s)-w(z,s))
 \mu_s(dw\times dz)\,ds
\\ \dsp \phantom{}
\ \ t\ge t_0\,,\ B\in{\cal B}(W).
\end{array} \]
This implies equality as a measure:
\begin{equation}
\eqna{Burgersprf16}
\begin{array}{l}\dsp
\mu_t(dw\times[y_C(\gamma,t),1))
=e^{-\tilde{\Omega}_w(t_0,t)}\mu_{t_0}(dw\times[y_0,1))
\\ \dsp \phantom{\biggl|}
+\int_{t_0}^t 
e^{-\tilde{\Omega}_w(s,t)}\int_{z\in[y_C(\gamma,s),1)} (w(1,s)-w(z,s))
 \mu_s(dw\times dz)\,ds.
\end{array}
\end{equation}
Using
\[ w(1,s)-w(z,s)=\int_z^1 \pderiv{w}{z}(x,s)\,dx \]
with Fubini's Theorem, we arrive at \eqnu{Burgersprf10}.
\QED
\end{prf}
Let us return to the proof of uniqueness in \thmu{Burgers}, and
suppose there is another pair $(\tilde{\mu},\tilde{y}_C)$
which satisfies all the properties stated in \thmu{Burgers}.
\lemu{Burgersprf} implies that $(\tilde{\mu}_t,\tilde{y}_C)$ satisfies
an integral equation similar to \eqnu{Burgersprf10},
\begin{equation}
\eqna{Burgersprf13}
\begin{array}{l}\dsp
\tilde{\mu}_t(dw\times[\tilde{y}_C(\gamma,t),1))
=e^{-\tilde{\Omega}_w(t_0,t)}\mu_{t_0}(dw\times[y_0,1))
\\ \dsp \phantom{\biggl|}
+\int_{t_0}^t 
e^{-\tilde{\Omega}_w(s,t)}\int_{x\in[\tilde{y}_C(\gamma,s),1)} \pderiv{w}{z}(x,s)
 \tilde{\mu}_s(dw\times [\tilde{y}_C(\gamma,s),x))\,dx\,ds.
\end{array}
\end{equation}
The first term in the right hand side is equal to that of
\eqnu{Burgersprf10}, because of the initial and boundary conditions
\eqnu{conservation} and \eqnu{soliditycond}.

Put
\begin{equation}
\eqna{Burgersprf11}
I(t)=\sup_{h}\sup_{y\in[0,1)} \biggl|\int_W h(w)
\tilde{\mu}_t(dw\times [y,1))-\int_W h(w)\mu_t(dw\times [y,1))\biggr|,
\end{equation}
and
\begin{equation}
\eqna{Burgersprf12}
J(t)=\sup_{h}\sup_{\gamma\in\Gamma_t} \biggl|\int_W h(w)
\tilde{\mu}_t(dw\times [\tilde{y}_C(\gamma,t),1))
-\int_W h(w)\mu_t(dw\times [y_C(\gamma,t),1))\biggr|,
\end{equation}
where the supremum for $h$ in the right hand sides
of \eqnu{Burgersprf11} and  \eqnu{Burgersprf12} are taken over
measurable functions $h:\ W\to\reals$ satisfying \eqnu{Burgersprf3}.
In particular, \eqnu{varphiyC} implies
\begin{equation}
\eqna{Burgersprf0}
\sup_{\gamma\in \Gamma_t}
|\tilde{y}_C(\gamma,t)-y_C(\gamma,t)|\le J(t).
\end{equation}

Fix $s$ and $x$, and put
\[
h_1(w)=\frac1{C_W}\,h(w)\,e^{-\tilde{\Omega}_w(s,t)}\,\pderiv{w}{z}(x,s).
\]
Then $\tilde{\Omega}_w(s,t)\ge0$ and \eqnu{Cw} imply that
$h_1:\ W\to\reals$ satisfies \eqnu{Burgersprf3} with $h=h_1$\,.
The definitions \eqnu{Burgersprf11} and \eqnu{Burgersprf12} then imply
\begin{equation}
\eqna{Burgersprf14}
\begin{array}{l}\dsp
\biggl|\int_W h(w)\,e^{-\tilde{\Omega}_w(s,t)}\,\pderiv{w}{z}(x,s)\,
\tilde{\mu}_s(dw\times[\tilde{y}_C(\gamma,s),x))
\\ \dsp \phantom{\biggl|}
-\int_W h(w)\,e^{-\tilde{\Omega}_w(s,t)}\,\pderiv{w}{z}(x,s)\,
\mu_s(dw\times[y_C(\gamma,s),x))\biggr|
\\ \dsp \phantom{}
\le
C_W\,\biggl|\int_W h_1(w)
\tilde{\mu}_s(dw\times[\tilde{y}_C(\gamma,s),1))
-\int_W h_1(w)\mu_s(dw\times[y_C(\gamma,s),1))\biggr|
\\ \dsp \phantom{\le}
+C_W\,\biggl|\int_W h_1(w)\tilde{\mu}_s(dw\times[x,1))
-\int_W h_1(w)\mu_s(dw\times[x,1))\biggr|
\\ \dsp \phantom{\biggl|}
\le C_W\,(I(s)+J(s)).
\end{array}
\end{equation}
Also, monotonicity of measure implies
\begin{equation}
\eqna{Burgersprf15}
\begin{array}{l}\dsp
\biggl|\int_W h(w)\,e^{-\tilde{\Omega}_w(s,t)}\pderiv{w}{z}(x,s)
 \tilde{\mu}_s(dw\times [\tilde{y}_C(\gamma,s),x))\biggr|
\\ \dsp {}
\le C_W \tilde{\mu}_s(W\times[0,1])=C_W\,.
\end{array}
\end{equation}
Substituting \eqnu{Burgersprf10} and \eqnu{Burgersprf13}
in \eqnu{Burgersprf12}, using \eqnu{Burgersprf14} and
\eqnu{Burgersprf15}, and also \eqnu{Burgersprf0},  we have
\begin{equation}
\eqna{Burgersprf17}
\begin{array}{l}\dsp
J(t)\le \sup_{h}\sup_{\gamma\in\Gamma_t}\int_{t_0}^t\biggl(
\biggl|\int_{\tilde{y}_C(\gamma,s)}^{y_C(\gamma,s)}C_W\,dx\biggr|
+ \int_{x\in[y_C(\gamma,s),1)}\biggl| C_W\,(I(s)+J(s))\biggr|\,dx\biggr)
\,ds
\\ \dsp \phantom{J(t)}
\le C_W\,\sup_{h}\sup_{\gamma\in\Gamma_t} \int_{t_0}^t
(|\tilde{y}_C(\gamma,s)-y_C(\gamma,s)|+I(s)+J(s))\,ds
\\ \dsp \phantom{J(t)}
\le  \int_0^t (I(s)+2J(s))\,ds.
\end{array}
\end{equation}

Next, since by assumption, for each $t$, 
the map $\gamma\mapsto y_C(\gamma,t)$ is surjective, 
we have, from \eqnu{Burgersprf11},
\begin{equation}
\eqna{Burgersprf18}
I(t)=\sup_{h}\sup_{\gamma\in\Gamma_t} \biggl|\int_W h(w)
\tilde{\mu}_t(dw\times [y_C(\gamma,t),1))
-\int_W h(w)\mu_t(dw\times [y_C(\gamma,t),1))\biggr|.
\end{equation}
Using Lipschitz continuity \eqnu{globalLipschitz},
the definition \eqnu{Burgersprf12}, and
\eqnu{Burgersprf0}, we have
\begin{equation}
\eqna{Burgersprf19}
I(t)\le \sup_{h}\sup_{\gamma\in\Gamma_t}
(|\tilde{y}_C(\gamma,t)-y_C(\gamma,t)|+J(t))
\le 2J(t).
\end{equation}
The inequalities \eqnu{Burgersprf18} and \eqnu{Burgersprf19}
imply
\[ J(t)\le 4\int_0^t J(s)\,ds,\ t\in[0,T], \]
which, by Gronwall's inequality, further implies $J(t)=0$, $t\in[0,T]$.
This proves $\tilde{y}_C=y_C$, and also \eqnu{Burgersprf19} now implies
$I(t)=0$, $t\in[0,T]$, which proves $\tilde{\mu}_t=\mu_t$\,.
This completes a proof of the uniqueness claim in \thmu{Burgers}.

\appendix

\section{Application of Schauder's fixed point theorem.}
\seca{Schauder}

In this section we consider the case where we keep the fundamental condition
\eqnu{rw}, but replace the global Lipschitz type condition \eqnu{Cw} by
a global bound condition on oscillation:
\begin{equation}
\eqna{Osc}
C'_W:= \sup_{w\in W} \sup_{(y,t),\,(y',t')\in[0,1]\times[0,T]}
|w(y,t)-w(y',t')|<\infty,
\end{equation}
and consider the existence of fixed points to the map
$G:\ \Theta_T\to\Theta_T$ defined by the explicit formula
\eqnu{HDLop_on_class_of_flows_f} and \eqnu{HDLop_on_class_of_flows_g},
where the notations are introduced in
\eqnu{Burgersassumps},
\eqnu{Burgersdensity1},
\eqnu{Gamma},
\eqnu{Gammat},
\eqnu{domyC},
\eqnu{class_of_conti_flows}, and
\eqnu{omegaOmega}.
Note that $G(\theta)\in\Theta_T$ holds with the weaker assumption
\eqnu{Osc}. 
This can be shown directly from the explicit expression
\eqnu{HDLop_on_class_of_flows_g}.
The only condition for $\Theta_T$ perhaps not 
obvious from the expression is the range condition
$G(\theta)(\gamma,t)\in [0,1]$, which can be shown as follows.
For $k=1,2,\ldots$ put
\[ \begin{array}{l}\dsp
I_k=
 \int_{0\le u_1\le \cdots\le u_k\le t_0}
 w(\theta((z,0),u_1),u_1)\, e^{-\Omega_{\theta,w,z}(0,u_1)}\,
\\ \dsp \phantom{I_k= \int_{0\le u_1\le \cdots\le u_k\le t_0}}
\times\prod_{i=2}^k \biggl( w(\theta((0,u_{i-1}),u_i),u_i)
\,e^{-\Omega_{\theta,w}(u_{i-1},u_i)}\biggr)\,
\\ \dsp \phantom{I_k= \int_{0\le u_1\le \cdots\le u_k\le t_0}}
\times e^{-\Omega_{\theta,w}(u_k,t)}
\prod_{i=1}^kdu_i
\end{array} \]
and 
\[ \begin{array}{l}\dsp
J_k=
 \int_{0\le u_1\le \cdots\le u_k\le t_0}
 w(\theta((z,0),u_1),u_1)\, e^{-\Omega_{\theta,w,z}(0,u_1)}\,
\\ \dsp \phantom{J_k= \int_{0\le u_1\le \cdots\le u_k\le t_0}}
\times\prod_{i=2}^k \biggl( w(\theta((0,u_{i-1}),u_i),u_i)
\,e^{-\Omega_{\theta,w}(u_{i-1},u_i)}\biggr)\,
\prod_{i=1}^kdu_i \,.
\end{array} \]
Note that non-negativity of $\Omega_{\theta,w}$ implies
$I_k\le J_k$.
We can perform the $u_k$ integration in $J_k$ to find
$J_k=J_{k-1}-I_{k-1}$,
which we can iterate to find
\[
\sum_{i=1}^k I_i= J_1-J_k+I_k \le J_1 = 1-e^{-\Omega_{\theta,w,z}(0,t_0)}.
\]
Substituting this in
\eqnu{HDLop_on_class_of_flows_g},
and using monotonicity of $\Omega_{\theta,w,z}$,
\eqnu{Burgersdensity1}, and $\lambda(W)=1$, we have
$1\ge G(\theta)((0,t_0),t)\ge 0$.
A similar estimate for \eqnu{HDLop_on_class_of_flows_f} is straightforward,
hence we conclude $G(\theta)\in \Theta_T$.
\begin{thm}
\thma{FPexistence}.
Under the condition \eqnu{rw} and \eqnu{Osc},
the map $G:\ \Theta_T\to\Theta_T$ has a fixed point.
\DDD
\end{thm}
\begin{rem}
Since the proof relies on Schauder's fixed point theorem, 
our proof has no control of uniqueness of fixed points.
\DDD\end{rem}

The map $G$ maps $\Theta_T$ into itself, and $\Theta_T$ is a
subset of a Banach space (with the supremum norm)
of continuous functions $C^0(\Delta_T;[0,1])$
taking values in a finite interval $[0,1]\subset\reals$.
The domain $\Delta_T$ is homeomorphic to a rectangle,
since its parameterization in the definition \eqnu{domyC} is 
homeomorphic to a trapezoid.

The Schauder fixed point theorem states \cite[(2.4.3)]{Berger} that
a compact map of a closed bounded convex set in a Banach space into itself
has a fixed point.
(The notational correspondence between here and \cite[\S 2.4]{Berger} is
given by $\dsp X=C^0(\Delta_T;\;[0,1])$, $K=U=\Theta_T$, and $f=G$.)
We have shown $G(\Theta_T)\subset \Theta_T$ at the beginning of
this section.

Concerning the required properties for the domain $\Theta_T$ of the map $G$,
we have noted that $C^0(\Delta_T;[0,1])$ is a bounded set.
For a sequence of continuous and monotone functions,
the limit function with respect to the supremum norm
also is continuous and monotone,
and since for $\theta\in\Theta_T$, $\theta((0,t),t)=0$ and
$\theta((1,0),t)=1$ holds, these properties are also preserved in the limit.
This and continuity imply surjectivity of the limit function.
Therefore $\Theta_T$ is a closed set.
The continuity, monotonicity, the properties
$\theta((0,t),t)=0$ and $\theta((1,0),t)=1$ are also preserved by
convex linear combination, hence $\Theta_T$ is also convex.
Thus $\Theta_T$ is a closed, bounded, convex set.

It remains to prove compactness of $G$.
Since $C^0(\Delta_T;\;[0,1])$ is a bounded set with respect to the supremum
norm, the Arzela-Ascoli theorem implies that
it is sufficient to prove (i) that the map $G:\ \Theta_T\to\Theta_T$
is continuous, and (ii) that the functions in the image set $G(\Theta_T)$ 
are equicontinuous, which we prove in \lemu{Gconti} and
\lemu{GTheta_is_equiconti}, respectively.

Note first that non-negativity of $w\in W$ obviously implies
\begin{equation}
\eqna{Gcontiprf0}
0<e^{-\tilde{\Omega}_w(s,t)}\le 1,\ s\le t,
\end{equation}
where $\tilde{\Omega}_w(s,t)$ is as in \eqnu{wstar}.
\begin{prp}
\prpa{Gcontiprf1}
For $\theta$ and $\theta'$ in $\Theta_T$\,, and
$\dsp (\gamma,t)\in\Delta_T$ with $\gamma=(z,s)$, we have
\begin{enumerate}
\item
$\dsp w(1,t)-C'_W\le w(\theta(\gamma,t),t) 
\le w(1,t)+C'_W\le \wnorm{w}+C'_W$,

\item
$\dsp 0<e^{-\Omega_{\theta,w,z}(s,t)}
\le e^{-\tilde{\Omega}_w(s,t)+C'_W\,(t-s)}$,

\item
$\dsp |e^{-\Omega_{\theta',w,z}(s,t)}-e^{-\Omega_{\theta,w,z}(s,t)}|
\le e^{-\tilde{\Omega}_w(s,t)+C'_W\,(t-s)}
\int_s^t |w(\theta'(\gamma,u),u)-w(\theta(\gamma,u),u)|\,du$,

\end{enumerate}
where $C'_W$ is as in \eqnu{Osc}, and
$\wnorm{w}$ is defined by \eqnu{Wnorm}.
\DDD
\end{prp}
\begin{prf}
\eqnu{Osc} implies
\[
|w(y,t)-w(1,t)|\le C'_W,\ w\in W,\ (y,t)\in[0,1]\times[0,T],
\]
which further implies
\[
 0<e^{-\Omega_{\theta,w,z}(s,t)}\le e^{-\int_s^t w(1,u)\,du+C'_W\,(t-s)},
\]
These estimates imply the first $2$ estimates.
The last estimate follows from these estimates and
\[
|e^{-x'}-e^{-x}|=|e^{-(x'\vee x)}-e^{-(x'\wedge x)}|
= e^{-(x'\wedge x)} |\, (1- e^{-|x'-x|})
\\ \dsp \le e^{-(x'\wedge x)} |x'-x|.
\]
\QED
\end{prf}

\begin{lem}
\lema{Gconti}
$G:\ \Theta_T\to\Theta_T$ is a continuous map.
\DDD
\end{lem}
\begin{prf}
Let $\dsp\theta,\theta'\in \Theta_T$, and put
$(\gamma,t)\in\Delta_T$ and $\gamma=(y_0,t_0)$.

If $\gamma=(y_0,0)\in\Gamma_i$ ($t_0=0$),
\eqnu{HDLop_on_class_of_flows_f},
\prpu{Gcontiprf1}C\eqnu{Gcontiprf0}, and \eqnu{Burgersdensity1}
imply
\begin{equation}
\eqna{Gcontiprf3}
\begin{array}{l}\dsp
\sup_{y_0\in[0,1]}\sup_{t\in[0,T]}
|G(\theta')((y_0,0),t)-G(\theta)((y_0,0),t)|
\\ \dsp {}
\le 
\sup_{y_0\in[0,1]}\sup_{t\in[0,T]}
\int_{W\times[y_0,1)} 
 e^{-\tilde{\Omega}_w(0,t)+C'_W\,t}
\int_0^t |w(\theta'((z,0),u),u)-w(\theta((z,0),u),u)|\,du
\\ \dsp \phantom{\sup_{y_0\in[0,1]}\sup_{t\in[0,T]}\int_{W\times[y_0,1)}}
\times\sigma(w,z)\,\lambda(dw)\,dz
\\ \dsp {}
\le  T\,e^{C'_W\,T}\int_{W\times[0,1)} 
\sup_{u\in[0,T]}\sup_{z\in[0,1]}
 |w(\theta'((z,0),u),u)-w(\theta((z,0),u),u)|
\,\lambda(dw).
\end{array}
\end{equation}
Concerning the rightmost hand side, we have
\[
\sup_{u\in[0,T]}\sup_{z\in[0,1]}
 |w(\theta'((z,0),u),u)-w(\theta((z,0),u),u)|
\le 2\wnorm{w},
\]
while \eqnu{rw} implies
$\dsp
\int_{W\times[0,1)} \wnorm{w}\lambda(dw)=M_W<\infty.
$
Hence the integrand in the right hand side of \eqnu{Gcontiprf3}
is bounded, pointwise in $w\in W$, uniformly in $\theta'$ by an integrable 
function. Therefore, thanks to dominated convergence theorem
we may interchange the order of integration and the limit
$\theta'\to\theta$ in the right hand side of \eqnu{Gcontiprf3}.

$W\subset C^1([0,1]\times[0,T];[0,\infty))$ and 
$[0,1]\times[0,T]$ is compact, hence each $w\in W$ is uniformly continuous.
Hence for any $\eps>0$, there exists $\delta>0$ such that
\[ (\forall y,y'\in[0,1];\ |y'-y|<\delta)(\forall u\in[0,T])
\ |w(y,u)-w(y',u)|<\eps. \]
If the supremum norm of $\theta'-\theta$ is less than $\delta$
we have
\[
\sup_{u\in[0,T]}\sup_{z\in[0,1]}
 |\theta'((z,0),u)-\theta((z,0),u)|<\delta,
\]
which further implies
\[
\sup_{u\in[0,T]}\sup_{z\in[0,1]}
 |w(\theta'((z,0),u),u)-w(\theta((z,0),u),u)|
\le\eps,
\]
hence \eqnu{Gcontiprf3} implies
\[
\limsup_{\theta'\to\theta}
\sup_{y_0\in[0,1]}\sup_{t\in[0,T]}
|G(\theta')((y_0,0),t)-G(\theta)((y_0,0),t)|
\le T\,e^{C'_W\,T}\eps.
\]
Since $\eps>0$ is arbitrary,
\begin{equation}
\eqna{Gcontiprf2}
\lim_{\theta'\to\theta}
\sup_{y_0\in[0,1]}\sup_{t\in[0,T]}
|G(\theta')((y_0,0),t)-G(\theta)((y_0,0),t)|=0.
\end{equation}

Next if $\gamma=(0,t_0)\in\Gamma_t\cap\Gamma_b$ ($y_0=0$),
we proceed as in exact analogy to the proof of \lemu{FPtheoremprf2},
to find 
\[\begin{array}{l}\dsp
\sup_{t\in[0,T]}\sup_{t_0\in[0,t]}
|G(\theta')((0,t_0),t)-G(\theta)((0,t_0),t)|
\\ \dsp {}
\le T\,e^{2C'_W\,T}
\int_W 
\sup_{(\gamma,u)\in\Delta_T}
|w(\theta'(\gamma,u),u)-w(\theta(\gamma,u),u)|
\,\lambda(dw).
\end{array}\]
Therefore,
as in the same reasoning as we derive \eqnu{Gcontiprf2} from
\eqnu{Gcontiprf3},
\begin{equation}
\eqna{Gcontiprf10}
\lim_{\theta'\to\theta}
\sup_{t\in[0,T]}\sup_{t_0\in[0,t]}
|G(\theta')((0,t_0),t)-G(\theta)((0,t_0),t)|=0.
\end{equation}

Finally, \eqnu{Gcontiprf2} and \eqnu{Gcontiprf10} imply
\[
\lim_{\theta'\to\theta}
\sup_{(\gamma,t)\in\Delta_T}
|G(\theta')(\gamma,t)-G(\theta)(\gamma,t)|=0,
\]
which proves the continuity of $G:\ \Theta\to\Theta$.
\QED
\end{prf}

\begin{lem}
\lema{GTheta_is_equiconti}
The functions in the set $G(\Theta_T)$ are equicontinuous.
\DDD
\end{lem}
\begin{prf}
We see from elementary calculus using the mean value theorem and
triangular inequality that the following uniform estimates
on the derivatives of $G(\theta)$ imply equicontinuity:
\begin{equation}
\eqna{GThetaequicontiprf1}
0\le \pderiv{}{y_0}G(\theta)((y_0,0),t)\le 1,
\ \ y_0\in[0,1],\ t\in[0,T],
\end{equation}
\begin{equation}
\eqna{GThetaequicontiprf2}
0\le \pderiv{}{t}G(\theta)((y_0,0),t)\le M_W+C'_W\,,
\ \ y_0\in[0,1],\ t\in[0,T],
\end{equation}
\begin{equation}
\eqna{GThetaequicontiprf3}
0\le -\pderiv{}{t_0}G(\theta)((0,t_0),t)\le (M_W+C'_W)\,e^{2C'_WT},
\ \ t_0\in[0,t],\ t\in[0,T],
\end{equation}
\begin{equation}
\eqna{GThetaequicontiprf4}
0\le \pderiv{}{t}G(\theta)((0,t_0),t)\le (M_W+C'_W)\,e^{2C'_WT},
\ \ t_0\in[0,t],\ t\in[0,T].
\end{equation}
The remainder of the proof is devoted to proving these estimates.

To prove \eqnu{GThetaequicontiprf1},
differentiate the explicit formula \eqnu{HDLop_on_class_of_flows_f}
by $t_0$ and use \eqnu{Gcontiprf0} and \eqnu{Burgersassumps}.
We have
\[ \begin{array}{l}\dsp
0\le \pderiv{}{y_0}G(\theta)((y_0,0),t)
=\int_{W}
e^{-\Omega_{\theta,w,y_0}(0,t)}\,\sigma(w,y_0)\,\lambda(dw)
\\ \dsp \phantom{0\le \pderiv{}{y_0}G(\theta)((y_0,0),t)}
\le \int_W \sigma(w,y_0)\lambda(dw) = 1,
\end{array} \]
which proves \eqnu{GThetaequicontiprf1}.

To prove \eqnu{GThetaequicontiprf2},
differentiate \eqnu{HDLop_on_class_of_flows_f} by $t$,
and use \eqnu{omegaOmega}, \prpu{Gcontiprf1}, \eqnu{Gcontiprf0},
\eqnu{Burgersdensity1},
and \eqnu{rw}, to find
\[ \begin{array}{l}\dsp
0\le \pderiv{}{t}G(\theta)((y_0,0),t)
=\int_{W\times[y_0,1)} w(\theta((z,0),t),t)\,
e^{-\Omega_{\theta,w,z}(0,t)}\,\sigma(w,z)\,\lambda(dw)\,dz
\\ \dsp \phantom{0\le \pderiv{}{t}G(\theta)((y_0,0),t)}
\le \int_{W\times[0,1)} (\wnorm{w}+C'_W)\, \sigma(w,z)\,\lambda(dw)\,dz
=M_W+C'_W\,,
\end{array} \]
which proves \eqnu{GThetaequicontiprf2}.

Proofs of \eqnu{GThetaequicontiprf3} and
\eqnu{GThetaequicontiprf4} are similar.
We differentiate \eqnu{HDLop_on_class_of_flows_g} by $t_0$ and $t$,
respectively, and follow a similar line. 
The only new point is that we apply
\eqnu{integralopcompactness} in a similar way as in the proof of
\lemu{FPtheoremprf2}.
This completes a proof of \lemu{GTheta_is_equiconti}.
\QED
\end{prf}

As discussed in the beginning of this appendix,
\lemu{Gconti} and \lemu{GTheta_is_equiconti} prove
\thmu{FPexistence}.

\end{document}